\documentclass[12pt,letterpaper]{article}

\usepackage{amsmath}
\usepackage[OT2,OT1]{fontenc}
\usepackage{epsfig,graphicx,epsf}
\usepackage{fancyhdr}
\usepackage{pstricks}
\usepackage{pst-plot}
\usepackage{amscd}
\usepackage{amssymb}
\usepackage{latexsym}

\newcounter{psfig}
\makeatletter
\def\cleardoublepage{\clearpage\if@twoside\ifodd\c@page\else
\hbox{}
\thispagestyle{empty}
\newpage
\if@twocolumn\hbox{}\newpage\fi\fi\fi}
\makeatother

\numberwithin{equation}{section}

\newcommand\gam{\gamma}         
\newcommand\del{\delta}         \newcommand\Del{\Delta}
\newcommand\eps{\varepsilon}

\newcommand\iot{\iota}
\newcommand\kap{\kappa}
\newcommand\lam{\lambda}        \newcommand\Lam{\Lambda}
\newcommand\sig{\sigma}         \newcommand\Sig{\Sigma}
\newcommand\vrh{\varrho}
\newcommand\vth{\vartheta}
\newcommand\vph{\varphi}
\newcommand\ome{\omega}

\newcommand\calC{{\mathcal{C}}}

\newcommand\calL{{\mathcal{L}}}

\newcommand\calN{{\mathcal{N}}}
\newcommand\calO{{\mathcal{O}}}

\newcommand\calV{{\mathcal{V}}}

\newcommand\calZ{{\mathcal{Z}}}

            \newcommand\rmH{{\mathrm H}}

\newcommand\RR{\mathbb{R}}
\newcommand\TT{\mathbb{T}}

\newcommand\ZZ{\mathbb{Z}}

\newcommand\CC{\mathbb{C}}

\newcommand\NN{\mathbb{N}}

\newcommand\grK{{\mathfrak{K}}} \newcommand\grk{{\mathfrak{k}}}

 \newcommand\grt{{\mathfrak{t}}}

\newcommand\ddd{,\dots,}
\renewcommand\d{\partial}

\newcommand\sdp{\times \hskip -0.48em {\raise 0.3ex
\hbox{$\scriptscriptstyle |$}}} % semidirect product

\newcommand\Aut{\operatorname{Aut}}

\newcommand\cd{\operatorname{cd}}

\newcommand\coker{\operatorname{coker}}

\newcommand\Crit{\operatorname{Crit}}
\newcommand\Der{\operatorname{Der}}

\newcommand\diag{\operatorname{diag}}

\newcommand\Diff{\operatorname{Diff}}

\newcommand\GL{\operatorname{GL}}

\newcommand\Id{\operatorname {Id}}
\newcommand\id{\operatorname {id}}
\newcommand\im{\operatorname {im}}
\newcommand\IM{\operatorname{Im}}

\newcommand\MOD{\,\operatorname{mod}\,}

\newcommand\MRep{M{\mathrm{-Rep}}}

\newcommand\rk{\operatorname{rk}}

\newcommand\SL{{\rm SL}}

\newcommand\SO{{\rm SO}}

\newcommand\Trace{\operatorname{Trace}}

\newcommand\oB{{\overline{B}}}

\renewcommand{\>}{\rangle}
\newcommand{\<}{\langle}
\newcommand{\half}{\frac{1}{2}}

\newtheorem{Thm}{Theorem}[section]
\newtheorem{Cor}[Thm]{Corollary}
\newtheorem{Lem}[Thm]{Lemma}
\newtheorem{Fac}[Thm]{Fact}
\newtheorem{Nt}[Thm]{Notation}
\newtheorem{Prop}[Thm]{Proposition}
\newtheorem{Conjec}[Thm]{Conjecture}

\newtheorem{Compt}[Thm]{Computation}
\newtheorem{Rem}[Thm]{Remark}

\newtheorem{Cex}[Thm]{Counterexample}
\newtheorem{Para}{\S}
\newtheorem{LemP}{Lemma}[section]
\newtheorem{PropP}[LemP]{Proposition}

\newtheorem{Def}[Thm]{Definition}

\newif\ifShowLabels
\ShowLabelstrue
\newdimen\theight
\def\TeXref#1{%
        \leavevmode\vadjust{\setbox0=\hbox{{\tt
                \quad\quad  {\small \textrm #1}}}%
        \theight=\ht0
        \advance\theight by \lineskip
        \kern -\theight \vbox to
        \theight{\rightline{\rlap{\box0}}%
        \vss}%
        }}%

\ShowLabelsfalse% comment this out if labels should be printed

\renewcommand{\sec}[2]{\section{#2}\label{S:#1}%
        \ifShowLabels \TeXref{{S:#1}} \fi}
\newcommand{\ssec}[2]{\subsection{#2}\label{SS:#1}%
        \ifShowLabels \TeXref{{SS:#1}} \fi}

\newcommand{\refs}[1]{Section ~\ref{S:#1}}
\newcommand{\refss}[1]{Paragraph ~\ref{SS:#1}}
\newcommand{\reft}[1]{Theorem ~\ref{T:#1}}

\newcommand{\refp}[1]{Proposition ~\ref{P:#1}}

\newcommand{\refd}[1]{Definition ~\ref{D:#1}}

\newcommand{\refas}[1]{Appendix ~\ref{S:#1}}

\newenvironment{thm}[1]%
        { \begin{Thm} \label{T:#1}  \ifShowLabels \TeXref{T:#1} \fi }%
        { \end{Thm} }

\renewcommand{\th}[1]{\setlength{\parindent}{0cm}\begin{thm}{#1} \it }
\renewcommand{\eth}{\setlength{\parindent}{0.8cm}\end{thm} }

\newenvironment{NT}[1]%
        { \begin{Nt} \label{N:#1}  \ifShowLabels \TeXref{N:#1} \fi }%
        { \end{Nt} }

\newcommand{\nt}[1]{\setlength{\parindent}{0cm}\begin{NT}{#1} \it }
\newcommand{\ent}{\setlength{\parindent}{0.8cm}\end{NT} }

\newenvironment{compt}[1]%
        { \begin{Compt} \label{CO:#1}  \ifShowLabels \TeXref{CO:#1} \fi }%
        { \end{Compt} }

\newcommand{\co}[1]{\setlength{\parindent}{0cm}\begin{compt}{#1} \it }
\newcommand{\eco}{\setlength{\parindent}{0.8cm}\end{compt} }

\newenvironment{FC}[1]%
        { \begin{Fac} \label{F:#1}  \ifShowLabels \TeXref{F:#1} \fi }%
        { \end{Fac} }

\newcommand{\fac}[1]{\setlength{\parindent}{0cm}\begin{FC}{#1} \it }
\newcommand{\efac}{\setlength{\parindent}{0.8cm}\end{FC} }

\newenvironment{CX}[1]%
        { \begin{Cex} \label{CX:#1}  \ifShowLabels \TeXref{CX:#1} \fi }%
        { \end{Cex} }

\newcommand{\cex}[1]{\setlength{\parindent}{0cm}\begin{CX}{#1} \it }
\newcommand{\ecex}{\setlength{\parindent}{0.8cm}\end{CX} }

\newenvironment{lemma}[1]%
        { \begin{Lem} \label{L:#1}  \ifShowLabels \TeXref{L:#1} \fi }%
        { \end{Lem} }
\newcommand{\lem}[1]{\setlength{\parindent}{0cm}\begin{lemma}{#1} \it}
\newcommand{\elem}{\end{lemma}\setlength{\parindent}{0.8cm}}

\newenvironment{lemmaP}[1]%
        { \begin{LemP} \label{LL:#1}  \ifShowLabels \TeXref{LL:#1} \fi }%
        { \end{LemP} }
\newcommand{\lemm}[1]{\setlength{\parindent}{0cm}\begin{lemmaP}{#1} \it}
\newcommand{\elemm}{\setlength{\parindent}{0.8cm}\end{lemmaP}}

\newenvironment{parag}[2]%
        { \begin{Para} \label{PA:#1}\bf{#2}\hspace{.3cm}
          \ifShowLabels \TeXref{PA:#1} \fi }%
        { \end{Para}\vspace{0.2cm} }
\newcommand{\para}[2]{\begin{parag}{#1}{#2} \rm}
\newcommand{\epara}{\end{parag}}

\newenvironment{propos}[1]%
        { \begin{Prop} \label{P:#1}  \ifShowLabels \TeXref{P:#1} \fi }%
        { \end{Prop} }
\newcommand{\prop}[1]{\setlength{\parindent}{0cm}\begin{propos}{#1}\it }
\newcommand{\eprop}{\end{propos}\setlength{\parindent}{0.8cm}}

\newenvironment{proposp}[1]%
        { \begin{PropP} \label{PP:#1}  \ifShowLabels \TeXref{PP:#1} \fi }%
        { \end{PropP} }
\newcommand{\propp}[1]{\setlength{\parindent}{0cm}\begin{proposp}{#1}\it }
\newcommand{\epropp}{\end{proposp}\setlength{\parindent}{0.8cm}}

\newenvironment{corol}[1]%
        { \begin{Cor} \label{C:#1}  \ifShowLabels \TeXref{C:#1} \fi }%
        { \end{Cor} }
\newcommand{\cor}[1]{\setlength{\parindent}{0cm}\begin{corol}{#1} \it }
\newcommand{\ecor}{\end{corol}\setlength{\parindent}{0.8cm}}

\newenvironment{conjec}[1]%
        { \begin{Conjec} \label{Conj:#1}  \ifShowLabels \TeXref{C:#1} \fi }%
        { \end{Conjec} }
\newcommand{\conj}[1]{\setlength{\parindent}{0cm}\begin{conjec}{#1} \it }
\newcommand{\econj}{\end{conjec}\setlength{\parindent}{0.8cm}}

\newenvironment{defeni}[1]%
        { \begin{Def} \label{D:#1}  \ifShowLabels \TeXref{D:#1} \fi }%
        { \end{Def} }
\newcommand{\defe}[1]{\setlength{\parindent}{0cm}\begin{defeni}{#1} \it }
\newcommand{\edefe}{\end{defeni}\setlength{\parindent}{0.8cm}}

\newenvironment{remark}[1]%
        { \begin{Rem} \label{R:#1}  \ifShowLabels \TeXref{#1} \fi \rm}%
        { \end{Rem} }
\newcommand{\rem}[1]{\setlength{\parindent}{0cm}\begin{remark}{#1}}
\newcommand{\erem}{\end{remark}\setlength{\parindent}{0.8cm}}

\newcommand{\Label}[1]%
        {\label{#1}%
            \ifShowLabels \TeXref{#1} \fi }%

\newcommand{\eq}[1]%
         {\newline \ifShowLabels \TeXref{{E:#1}} \fi
                \begin{equation} \label{E:#1}}
\newcommand{\eeq}{\end{equation}}

\newcommand{\eprf}{\nopagebreak\end{proof} }

\newcommand{\pf}{{\noindent\it Proof. }}

\newcommand{\epf}{\protect\nopagebreak\hspace*{0.2cm}\hfill$\Box$\\[0.5ex plus
  0.5ex minus 0.2ex]}

\newcommand\beqn{\begin{equation}}
\newcommand\eeqn{\end{equation}}

\newcommand{\bear}{\begin{eqnarray*}}
\newcommand{\eear}{\end{eqnarray*}}

\newcommand{\bEAR}{\begin{eqnarray}}
\newcommand{\eEAR}{\end{eqnarray}}

\newcommand{\nid}{\noindent}
\newcommand{\pid}{\hspace*{0.6cm}}

\newcommand{\lra}{\longrightarrow}

\begin{document}

\vspace*{1.5cm}
{\LARGE{
\centerline{$T^3$-fibrations on compact six-manifolds}
}}
\vspace{1cm}
{\large{\centerline{Patrick Baier}}}
\vspace{0.4cm}
{\large \it\centerline{Dept. of Mathematics, UC Riverside, CA 92521}}
\vspace{1cm}

% ^

{\bf\centerline{Abstract}}
\begin{center}
\begin{minipage}{11.8cm}
We describe a simple way of constructing torus fibrations $T^3\to X\to S^3$ 
which degenerate canonically
over a knot or link $\grk\subset S^3$. We show that the topological invariants
of $X$ can be computed algebraically from the monodromy representation of
$\pi_1(S^3\setminus\grk)$ on $H_1(T^3,\ZZ)$. 
We use this to obtain 
some previously unknown $T^3$-fibrations $S^3\times 
S^3\to S^3$ and $(S^3\times S^3)\#(S^3\times S^3)\#(S^4\times S^2)\to S^3$ 
whose discriminant locus is a torus knot $\grt(p,q)\subset S^3$. 
\end{minipage}
\end{center}
\vspace{1cm}

\sec{S1}{Introduction}

Our study of $T^3$-fibrations in dimension six is strongly motivated by
recent advances surrounding the SYZ conjecture in mirror symmetry, especially
Mark Gross' work on special Lagrangian fibrations \cite{G1,G3}. Gross defines
a class of $T^3$-fibrations $X\to S^3$ which he calls ``well-behaved 
and admissible''
and shows that, at a topological level, they display the properties the 
SYZ-conjecture leads us to expect
of special Lagrangian fibrations on Calabi-Yau manifolds. 
In this paper we investigate a special case
of Gross' well-behaved, admissible fibrations. They develop 
canonical singularities
over a smoothly embedded knot or link $\grk\subset
S^3$. However, it turns out that the total spaces of our $T^3$-fibrations 
are often rather far from being Calabi-Yau which places this work
outside the immediate context of mirror symmetry, as an independent 
piece of differential geometry. \\
\pid An {\em affine $T^n$-bundle} is a fibre bundle 
$f\colon X\to B$ whose fibres
are $n$-tori,
equipped with an affine structure. More precisely, if we write 
$\TT^n=\RR^n/\ZZ^n$ for the standard linear torus, an
affine structure on the $T^n$-bundle $X$ can be specified by 
giving local trivialisations 
\[\vph_i\colon X|_{U_i}\lra U_i\times\TT^n\]
for some cover $\{U_i\}$ of $B$ so that $\vph_j\circ\vph_i^{-1}$ restricts to
an affine transformation on each torus fibre $\{x\}\times\TT^n$ with $x\in 
U_i\cap U_j$. We always assume that the total space $X$ and 
the base $B$ are orientable and connected.
An affine $T^n$-bundle
$f\colon X\to B$ is determined, up to equivalence, by its monodromy 
representation 
\[\vrh\colon \pi_1(B,b)\to\SL(H_1(F_b,\ZZ))\cong\SL(n,\ZZ)\] 
\nid where $b\in B$ is any basepoint and $F_b=f^{-1}(b)$, together with 
its Chern class 
\[c(f)\in H^2(B,\calL)\]
\nid  where $\calL$ is a bundle of lattices with $\calL_b=H_1(F_b,\ZZ)
\cong\ZZ^n$ and monodromy $\vrh$.
The affine $T^n$-bundle $f\colon X\to B$ admits a section
if and only if $c(f)=0$. In this case it is 
called a {\em linear $T^n$-bundle}, 
where we think of the section as determining an 
origin in each fibre.
Such a linear $T^n$-bundle $X$ with section can be recovered 
from its monodromy 
representation $\vrh$ as the bundle 
\[X=\widehat{B}\times_\vrh\TT^n\]
\nid  where $\widehat{B}$ 
is the universal cover of $B$ and $\pi_1(B,b)$ acts on $\widehat{B}$ by deck
transformations and on $\TT^n$ via $\vrh$, preserving the linear structure 
on $\TT^n$. \\
\hspace*{4ex}Not every $T^n$-bundle 
with section can be equipped with such a linear structure on the fibres. For 
instance, there are no non-trivial monodromy representations 
and Chern-classes on $S^m$ for $m>2$, so every affine $T^n$-bundle 
$X\to S^m$ is trivial, $X\cong S^m\times\TT^n$. General $T^n$-bundles over 
$S^m$ on the other side 
are classified by $\pi_{m-1}(\Diff(T^n))$ which is known to be 
non-trivial for certain values of $m,n$.\\
\hspace*{4ex}Few manifolds support the
structure of a linear or affine 
$T^n$-bundle. We have a better chance of finding, on
a given manifold $X$, an affine
$T^n$-{\em fibration}. This is roughly speaking an affine $T^n$-bundle with 
some singular fibres. We define

\defe{D1}{An \underline{affine $T^n$-fibration}
is a proper smooth map 
$f\colon X\to B$ between smooth manifolds so that for some proper closed
submanifold $\grk\subset B$ the restriction $f|_{B\setminus \grk}$ is an 
affine $T^n$-bundle.}\edefe

The assumption that the discriminant locus $\grk$  (that is the
set of critical values of $f$) 
be a submanifold makes this a somewhat unorthodox definition of a fibration, 
but this paper is only concerned with fibrations in the sense of \refd{D1}, 
and we do not want to invent a new name for them.\\

In \refs{S2}  we show that a linear $T^3$-bundle
$f_0\colon X_0\to B_0=S^3\setminus \grk$ (where $\grk$ is a compact 
one-dimensional submanifold, in other words, a knot or link) whose 
monodromy representation 
satisfies some algebraic condition can be canonically compactified to a 
$T^3$-fibration $f\colon X\to B$.

In \refs{S3} and \refas{A2} we study the relevant aspects
of the representation theory of knot (link) groups, in particular, 
for torus knots
$\grt(p,q)$ we determine all representations $\vrh\colon G(p,q)\to \SL(n,\ZZ)$ 
of the corresponding knot
groups $G(p,q)=\pi_1(S^3\setminus\grt(p,q))$ whose associated
$T^3$-fibrations can be compactified in this way.

In \refs{S4} we show how to compute the topological
invariants of the total space $X$ of such a fibration 
by purely algebraic means from the monodromy 
representation $\vrh$.

In \refs{S5} we apply the techniques developed in \refs{S4} 
to work out some concrete examples.
We construct an infinite number of non-isomorphic 
affine $T^3$-fibrations on the 
manifolds $S^3\times S^3$ and 
$(S^3\times S^3)\#(S^3\times S^3)\#(S^4\times S^2)$. 
The critical locus is a torus knot 
$\grt(2p',3q')\subset S^3$ with $\gcd(2p',3q')=1$, where $p'\in\NN$ is odd 
in the former case 
and is even in the latter case.\\

\nid {\bf Acknowledgments.\;} 
 I would like to thank R Herrera, N Hitchin, 
 M Lackenby and I Smith for helpful conversations, and YS Poon for inviting 
 me to Riverside.

% ^ 

\sec{S2}{The construction}

We describe a construction of $T^3$-fibrations using the ideas of \cite{G3}. 
Let the matrix
\begin{equation}\label{eq1}
\tilde{A}=\left(\begin{array}{cc}1&1\\0&1\end{array}\right)
\end{equation}
be given and define an action of the infinite cyclic group 
$\<\tilde{A}\>\cong\ZZ$ generated by $\tilde{A}$ 
on $\CC\times\TT^2$ by the formula
\begin{equation}\label{eqq2}
\tilde{A}^k\cdot(z,\xi)= (z+2\pi ik,\tilde{A}^k\xi)
\end{equation}
where $\tilde{A}\in\SL(2,\ZZ)$ acts on $\TT^2=\RR^2/\ZZ^2$ in the 
usual way. The quotient
\begin{equation}\label{eq3}
S_0=(\CC\times\TT^2)/{\<\tilde{A}\>}
\end{equation}
is a linear $T^2$-bundle over $\CC^*$ 
with section. Its monodromy with respect to the obvious basis of 
$H_1(\TT^2,\ZZ)$ is
given by (\ref{eq1}). 
It is well-known that there is a 
unique compactification 
\[S_0\subset S\stackrel{p}{\to}\CC\] 
which is topologically 
equivalent to an $I_1$ 
singularity in 
Kodaira's classification of singularities in elliptic fibrations. 
The singular fibre $F_0$ over $0\in\CC$ is a $T^2$ with one circle collapsed 
to a point. We are 
going to construct $T^n$-fibrations $f\colon X^{n+k}\to B^k$ 
with the property that
every $b\in B^k$ for which $F_b=f^{-1}(b)$ is singular 
has a neighbourhood $b\in U\subset B^k$ such that
there is a commutative diagram
\begin{equation}\label{Diag}
\begin{CD}f^{-1}(U)@>\Phi>\cong>S\times (\TT^{n-2}\times\RR^{k-2})\\
@VV{f}V @VV{(p,\mathrm{pr}_2)}V\\
U @>\phi>\cong> \CC\times\RR^{k-2}.\end{CD}
\end{equation}

\nid The fibration (\ref{Diag}) has a section.
With view towards our application on six-manifolds we now 
restrict to $T^3$-fibrations over $S^3$, thus $n=k=3$. 
However, much of this has 
obvious generalizations to arbitrary base manifolds and higher dimensions.
If a fibration has singularities of type (\ref{Diag}), then 
the discriminant locus $\grk=\phi^{-1}(\{0\}\times\RR)$ 
(the set of critical values 
of $f$) is an embedded one-dimensional submanifold. Because $\grk\subset 
S^3$ must be 
closed and $S^3$ is compact, $\grk$ is compact as well. 
This motivates the following definition.

\defe{D3}{A \underline{link}
is a smoothly embedded compact one-dimensional submanifold $\grk\subset S^3$. 
A link with just one connected component is a \underline{knot}.}\edefe 

Knot theorists often prefer to work in the category of 
piecewise linear spaces and maps. Also we have chosen to define a knot as
a submanifold rather than an equivalence class of oriented submanifolds
under the action 
of $\Diff^+(S^3)$, however, all constructions in what follows will only 
depend on the equivalence class of the knot or link. The idea is to 
compactify a $T^3$-bundle $T^3\to X_0\to S^3\setminus\grk$ over the complement
of a knot or link by ``gluing in'' singular fibres of type (\ref{Diag}) along
$\grk$ so that the section and the linear structure are preserved. 
We need to explain under which circumstances this is possible.

 \defe{D2}{Suppose $\grk\subset S^3$ is a knot or link and $b\in S^3\setminus 
 \grk$. An element $m\in\pi_1(S^3\setminus\grk,b)$ is called a 
 \underline{meridian} if there is a small open ball 
 $U\subset S^3$ for which $U\cap\grk$ is a connected arc and if we can choose 
 an $x\in U\setminus\grk$, 
 a loop $\gam_1\colon[0,1]\to U\setminus\grk$ based at $x$ representing a 
 generator of $\pi_1(U\setminus\grk,x)\cong\ZZ$ and a 
 path $\gam_2\colon[0,1]\to 
 S^3\setminus\grk$
 with $\gam_2(0)=b,\;\gam_2(1)=x$ so that $m$ is represented by 
 $\gam_2^{-1}*\gam_1*\gam_2$. An element $\ell\in\pi_1(S^3\setminus\grk,b)$ 
 is called a \underline{longitude} if there is a (``Seifert''-) 
 surface $\Sig\subset S^3$ with
 $\d\Sig$ a connected component of $\grk$ and a loop 
 $\lam_1\colon [0,1]\to\Sig$ which is based at some 
 $x\in\Sig$ and which is homotopic inside $\Sig$ 
 to $\d\Sig=\grk$  
 so that $\ell$ is represented by 
 $\lam_2^{-1}*\lam_1*\lam_2$ for some path $\lam_2$ from $b$ to $x$.
}\edefe

\nid Using this language, we can see that if $f\colon X\to S^3$ is a
$T^3$-fibration whose singularities are locally 
equivalent to (\ref{Diag}), the 
monodromy around any meridian must necessarily be
conjugate to
\begin{equation}\label{eeaa}
A=\left(\begin{array}{ccc}1&1&0\\0&1&0\\0&0&1\end{array}\right)\in\SL(3,\ZZ).
\end{equation}

\th{tcf}{Suppose $X_0\to S^3\setminus 
\grk$ is a linear $T^3$-bundle with section
whose monodromy $\vrh$ sends all meridians into the conjugacy class of 
the matrix $A$ given in (\ref{eeaa}). Then there exists a compactification
$X\to S^3$ with singularities of the form (\ref{Diag}).}\eth

\pf Suppose we have chosen
$b\in S^3\setminus \grk$, a loop $\mu$ which is based at $b$ and 
represents a meridian $m=[\mu]\in\pi_1(S^3\setminus\grk,b)$,
and an isomorphism $\SL(H_1(F_b,\ZZ))\cong\SL(3,\ZZ)$
which identifies the monodromy representation with a homomorphism
\[\vrh\colon\pi_1(S^3\setminus\grk,b)\to\SL(3,\ZZ)\]
such that $\vrh([\mu])=A$. 
Then let $B=\vrh([\lam])$ be the image of the longitude $\ell=[\lam]$ under $\vrh$. 
Obviously $[\mu][\lam]=[\lam][\mu]$, because if $\grK$ is an open 
``thickening'' of the knot then $\d(S^3\setminus\grK)\cong T^2$ and the 
inclusion of the boundary
$\iot\colon T^2\hookrightarrow S^3\setminus\grK$
satisfies $\iot_*(\pi_1(T^2))=\<[\lam],[\mu]\>$, and hence $[\lam]$ and
$[\mu]$ commute because $\pi_1(T^2)$ is Abelian. Consequently 
\[[A,B]=[\vrh([\lam]),\vrh([\mu])]=\vrh([[\lam],[\mu]])=0.\]
Then from definition (\ref{eqq2}) we can see that 
the standard action of 
$B$ on $\TT^3=\TT^2\times\TT^1$ descends to the quotient (\ref{eq3}) 
and induces an isomorphism $\vph_B$ of the linear $T^3$-bundle
on $S_0\times\TT^1\times\RR\to\CC^*\times\RR$.

 \prop{P1}{The map $\vph_B$ extends to an isomorphism of the
 $T^3$-fibration $S\times\TT^{1}\times\RR\to 
 \CC\times\RR$ which preserves the section.}\eprop

This can be checked explicitly by writing down the action of 
$\vph_B$ in our local model. 
We can now define a $\ZZ$-action 
on the linear $T^3$-fibration $S\times(\TT^1\times \RR)$ 
by the formula
\[k\cdot(x,\vth,t)= \vph_B^k(x,\vth,t+k).\]

\nid By dividing out this $\ZZ$-action 
we obtain a $T^3$-fibration $X_1\to \CC\times S^1$ with section whose smooth 
part is a $T^3$-bundle over $\CC^*\times S^1$ with monodromy given, 
in a suitable basis, by 
$A,B$. By construction, the singularities are locally equivalent to 
(\ref{Diag}).
Since a $T^n$-bundle with section is defined by its monodromy 
representation we can now compactify $X_0$ by gluing in (some
open subset of) $X_1$ using an isomorphism of affine $T^3$-bundles 
which preserves 
the linear structure and section.\epf

\defe{D2222}{We call fibrations of this kind \underline{good}
$T^3$-fibrations.}\edefe

\nid Thus a good $T^3$-fibration always has a section and an affine 
structure on its fibres and hence is 
completely determined by its monodromy representation.

\sec{S3}{Representations of knot groups}

In the previous section we described a construction of $T^3$-fibrations 
with section
degenerating over a knot or link $\grk\subset S^3$, using a gluing 
construction to compactify
a (linear) $T^3$-bundle over the link complement whose monodromy around each 
meridian was given by (\ref{eeaa}). We define

\defe{D5}{Let $\grk\subset S^3$ be a knot or link, $b\in S^3\setminus\grk$.
   By an {\underline{$M$-homomorphism}}
   we mean a homomorphism $\varrho\colon
   G(\grk)\to \SL(3,\ZZ)$ of the group
   $G(\grk)=\pi_1(S^3\setminus\grk,b)$ of $\grk$ which has the
   property that
   each meridian is mapped to the conjugacy class
   $\calC(A)$ of the matrix $A$ in (\ref{eeaa}). 
   More generally, if $\Lam$ is a free $\ZZ$-module of
   rank $\rk(\Lam)=3$, we call $\vrh\colon G(\grk)\to\SL(\Lam)$ an
   $M$-homomorphism if $\psi\vrh\psi^{-1}$ is an $M$-homomorphism
   for some (hence any) trivialisation $\psi\colon \Lam\to \ZZ^3$.
   Two homomorphisms $\varrho,\varrho'\colon
   G(\grk)\to \SL(\Lam)$ are conjugate if there is a fixed $M\in\SL(\Lam)$
   with $\varrho'=M^{-1}\varrho M$.
   An {\underline{$M$-re\-pre\-sen\-ta\-tion}} is a conjugacy class
   of $M$-ho\-mo\-mor\-phisms. }\edefe

The only input data for our
construction of $T^3$-fibrations
is an $M$-representation $\vrh\colon G\to \SL(3,\ZZ)$.
We will now study these $M$-representations of knot groups.

\defe{D0}{We call an element $\xi=(\xi_1,\xi_2,\xi_3)^{T}\in\ZZ^3$ 
   {\underline{primitive}} if 
\[\gcd(\xi_1,\xi_2,\xi_3)=1,\] 
   that is, if
   it is not a non-trivial 
   multiple of any other element in the lattice.}\edefe

\nid Define the usual bilinear form on $\ZZ^3$ by
\[\<\eta,\xi\>=\eta_1\xi_1+\eta_2\xi_2+\eta_3\xi_3.\]

\nid The 
following proposition is an easy to prove characterization of the 
conjugacy class $\calC(A)\subset\SL(3,\ZZ)$. 

\prop{PP}{A matrix $\tilde{A}\in\SL(3,\ZZ)$ is conjugate to $A\in 
\SL(3,\ZZ)$ if and only if it is of the form
\begin{equation}\label{ephi}
\tilde{A}=\Id+\eta\xi^T
\end{equation}

where $\eta,\xi\in\ZZ^3$ are primitive column vectors and $\<\eta,\xi\>=0$.
The pair $(\xi,\eta)$ is uniquely determined by $\tilde{A}$ up to
a common sign.}\eprop

\pf If $\tilde{A}$ and $A$ are conjugate then 
\[\rk(\tilde{A}-\Id)=\rk(A-\Id)=1\]

\nid and thus $\tilde{A}=\Id+\eta\xi^T$ for some $\eta,\xi\in\ZZ^3$. Moreover
\[3+\<\eta,\xi\>=\Trace(\tilde{A})=\Trace(A)=3\]

\nid and so $\<\eta,\xi\>=0$ as desired. Finally since 
$\IM(A-\Id)\cong\ZZ e_1$ 
is generated by a primitive element so is $\IM(\tilde{A}-\Id)$ 
which implies that
$\eta,\xi$ must be primitive.\\
\hspace*{4ex}Conversely let $\tilde{A}=\Id+\eta\xi^T$ 
be given with 
these properties.
We have to find $M\in\SL(3,\ZZ)$ such that
\[M^{-1}\eta\xi^{T}M=A-\Id=e_1e_2^{T}.\]

\nid This amounts to solving the equations 
\[M^{-1}\eta=e_1,\;\;\;M^{T}\xi=e_2.\]

\nid
The first of these equations simply says that $\eta$ has to be the
first column in $M$. Since $\xi$ is primitive, we can find
$\rho=(\rho_1,\rho_2,\rho_3)^T\in\ZZ^3$ such that
\[\langle\rho,\xi\rangle=\rho_1\xi_1+\rho_2\xi_2+\rho_3\xi_3=1.\]

\nid Let $\kappa=\eta\times\rho$ where the cross product is as in
$\RR^3$. Then $\kappa$ is primitive since if $\kappa$ was divisible,
then so would be $\kappa\times\xi$, however, 
\[\kappa\times\xi=(\eta\times\rho)\times\xi=\langle\eta,\xi\rangle\rho -
\langle\rho,\xi\rangle\eta=-\eta\]

\nid is primitive where we have used that $\langle\eta,\xi\rangle=0$. 
Hence we can
solve the equation 
$\langle\tilde{\sig},\kappa\rangle=1$ 
by some $\tilde{\sig}\in\ZZ^3$. But the latter sum is just
$\det(\eta,\rho,\tilde{\sig})=1$. Now set
$m=\langle\tilde{\sig},\xi\rangle$. Then if we replace
$\sig=\tilde{\sig}-m\rho$ and set 
\[M=(\eta,\rho,\sig)=\left(\begin{array}{ccc}\eta_1&\rho_1&\sig_1\\
\eta_2&\rho_2&\sig_2\\
\eta_3&\rho_3&\sig_3\end{array}\right)\]

\nid we will have
\[\det(M)=\det(\eta,\rho,\sig)=\det(\eta,\rho,\tilde{\sig}-m\rho)
=\det(\eta,\rho,\tilde{\sig})=1\]

\nid and the desired equations $Me_1=\eta$ and $M^{T}\xi=e_2$
follow. It is obvious that $\xi,\eta$ are uniquely determined up to a 
common sign.\epf

\prop{P23459}{Take $(\eta,\xi),(\tilde{\eta},{\tilde{\xi}})$
  with the properties given in \refp{PP}. Let
  $\tilde{A}=\Id+\tilde{\eta}{\tilde{\xi}}^T$ and $M=\Id+\eta\xi^T$. Then 
\[M^{-1}\tilde{A}M=\Id+\big(\tilde{\eta}-\langle\xi,\tilde{\eta}\rangle\eta 
\big)
\big({\tilde{\xi}}+\langle{\tilde{\xi}},\eta\rangle\xi\big)^T.\]
}\eprop

\pf The proof is straightforward. Note that $M^{-1}=\Id-\eta\xi^T$. 
\begin{eqnarray*}
&&M^{-1}\tilde{A}M\\
&=&(\Id-\eta\xi^{T})(\Id+\tilde{\eta}{\tilde{\xi}}^{T})
    (\Id+\eta\xi^{T})\\  
&=&\Id+\tilde{\eta}{\tilde{\xi}}^{T}-(\eta\xi^{T}) 
    (\tilde{\eta}{\tilde{\xi}}^{T})
    +(\tilde{\eta}{\tilde{\xi}}^{T})(\eta\xi^{T})
    -(\eta\xi^{T})(\eta\xi^{T})
    -(\eta\xi^{T})(\tilde{\eta}
    {{\tilde{\xi}}}^{T})(\eta\xi^{T})\\ 
&=&\Id+\big[\tilde{\eta}{\tilde{\xi}}^{T}
  -\langle\xi,\tilde{\eta}\rangle\eta{\tilde{\xi}}^{T}
  +\langle{\tilde{\xi}},\eta\rangle\tilde{\eta}\xi^{T}
  -\langle\xi,\tilde{\eta}\rangle \langle{\tilde{\xi}},\eta\rangle
  \eta\xi^{T}
  \big]\\
&=&\Id+\big(\tilde{\eta}-\langle\xi,\tilde{\eta}\rangle\eta\big)
   \big({\tilde{\xi}}+\langle{\tilde{\xi}},\eta\rangle\xi\big)^{T}.
\end{eqnarray*}

\nid The expression is invariant under the substitutions
$(\xi,\eta)\to (-\xi,-\eta)$ and $(\tilde{\xi},\tilde{\eta})\to
(-\tilde{\xi},-\tilde{\eta})$ and hence well-defined. \epf

We can use this formalism to find the $M$-representations of knot groups, 
starting from any Wirtinger presentation. A Wirtinger presentation describes
a knot or link group $G$ in terms of a set of meridians $g_1\ddd g_n$, 
subject to relations which can be read off a link diagram.
The meridians $g_1\ddd g_n$
all have to be mapped into the conjugacy class of the matrix $A$ in 
(\ref{eeaa}) under an $M$-representation $\vrh$. So we can make an Ansatz 
\[ \vrh(g_i)=\Id + \xi_i\eta_i^T,\;\;\;i=1\ddd n, \]

\nid where the $(\xi_i,\eta_i)$ are required to be pairs of primitive
orthogonal elements of $\ZZ^3$.
We translate the relations given by a Wirtinger presentation of $G$ into
a set of equations for the pairs $(\eta_i,\xi_i)$. The procedure is best 
explained by an example.

\th{P2345610}{Let $\grk=\grt(2,3)\subset S^3$ be the trefoil knot,
  $b\in S^3\setminus\grk$  and
  $G(\grk)=\pi_1(S^3\setminus \grk,b)$. Then
  the $M$-re\-pre\-sen\-ta\-tions of $G(\grk)$ can be classified in the
  following way:
There is one trivial $M$-re\-pre\-sen\-ta\-tion which
factorizes through the Abelianisation
  $G(\grk)^{\mathrm{ab}}=G(\grk)/[G(\grk),G(\grk)]\cong\ZZ$, \label{abel}
and one non-trivial $M$-re\-pre\-sen\-ta\-tion whose image in $\SL(3,\ZZ)$
is isomorphic to $\SL(2,\ZZ)$ (see (\ref{evrh})).}\eth

\pf As we can read from Figure 1, $G(\grk)$  has a 
Wirtinger presentation
\[G(\grk)=\big\langle\, g_1,g_2,g_3\;\big|\;
  g_2g_1=g_1g_3=g_3g_2\,\big\rangle.\]

\nid We rewrite these relations as
\begin{eqnarray}
\label{EW1}g_1^{-1}g_2g_1&=&g_3\\
\label{EW2}g_3^{-1}g_1g_3&=&g_2\\
\label{EW3}g_2^{-1}g_3g_2&=&g_1.
\end{eqnarray}

\begin{pspicture}(1,0)(16,4)
\pscustom[linewidth=1.5pt]{\pscurve(8.3,3.5)(8.7,3.6)(9.4,
  3.4)(9.9,2.7)(9.8,1.9)(9.4,1.3)(8.9,1.1)(8.4,1.0)(7.9,1.0)(7.5,1.0)
  (7.0,1.1)}
\pscustom[linewidth=1.5pt]{\pscurve(6.4,1.5)
  (6.1,2.2)(6.2,2.9)(6.8,3.5)(7.5,3.6)(8.2,3.4)(8.6,3.0)(8.8,2.6)  
  (9.0,2.3)(9.2,1.8)(9.2,1.5)}
\pscustom[linewidth=1.5pt]{\pscurve(9.2,0.9)(9.1,0.5)(8.4,0.1)(7.7,0.0)
  (7.0,0.5)(6.7, 1.1)  
  (6.8,1.8)(7.0,2.2)(7.2,2.6)(7.4,2.9)(7.6, 3.1)(7.8,3.3)}
\psline[linewidth=0.8pt,arrows=->](8.0,2.0)(7.4,2.4)
\pscurve[linewidth=0.8pt](7.0,2.6)(6.9,2.7)(6.85,2.8)(8.0,2.0)
\psline[linewidth=0.8pt,arrows=->](8.0,2.0)(8.6,2.4)
\pscurve[linewidth=0.8pt](9.0,2.6)(9.1,2.7)(9.15,2.8)(8.0,2.0)
\psline[linewidth=0.8pt,arrows=->](8.0,2.0)(8.1,1.3)
\pscurve[linewidth=0.8pt](8.15,0.9)(8.0,0.7)(7.85,0.8)(8.0,2.0)
\rput(6.5,2.8){$g_1$}
\rput(9.1,3.1){$g_2$}
\rput(8.3,0.6){$g_3$}
\rput(8.3,1.8){$b$}
\rput(8.0,-0.4){{\bf Figure 1: }
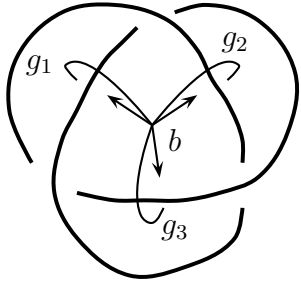
{\it The trefoil knot $\grt{(2,3)}$ (or $3_1$) and its
    meridians} $g_1,g_2,g_3$.} \label{pic}
\end{pspicture}
\vspace*{0.8cm}

\nid Note that $G(\grk)$ is in fact generated by $g_1$ and $g_2$ alone due
to (\ref{EW1}). Also the relation (\ref{EW3}) is a consequence of
(\ref{EW1}) and (\ref{EW2}) so that we only need to satisfy the first
two. We use the convention that for a manifold
with basepoint $(B,b)$ and elements $\gam_1,\gam_2$\label{loop} of its 
loop-group
the composition $\gam_1*\gam_2$ is the loop going {\em first} through
$\gam_2$ and {\em then} through $\gam_1$. This way our 
$M$-re\-pre\-sen\-ta\-tions $\pi_1(S^3\setminus\grk,b)\to\SL(3,\ZZ)$ will be
homomorphisms rather than anti-homomorphisms if we let $\SL(3,\ZZ)$ act
on $\ZZ^3$ from the {\em left}, as usual.\\

\pid We first show that for $\Lam\cong\ZZ^3$, every non-trivial 
$M$-homomorphism on $\Lam$ belongs to a family $\varrho_{abc}\colon
G(\grk)\to \SL(\Lam)$ which is parametrised bijectively by all 
$(a,b,c)\in\ZZ^3$. 
Given any $M$-homomorphism $\vrh_\Lam\colon G(\grk)\to\SL(\Lam)$ we can find
a trivialisation
$\Lam\to\ZZ^3$ so that $\vrh_\Lam$ becomes an $M$-homomorphism
$\vrh\colon G(\grk)\to\SL(3,\ZZ)$ of the form (see (\ref{ephi}))
\begin{eqnarray}
\label{eqm}  \varrho(g_1)&=&\Id+e_1e_2^T\,\;=\;\,A\\
\label{eqmm} \varrho(g_2)&=&\Id+\eta\xi^T
\end{eqnarray}
\nid for some $(\eta,\xi)$ primitive and orthogonal. Then $\varrho(g_3)$ 
is determined by
(\ref{EW1}). Namely, using \refp{P23459} and (\ref{eqm}), (\ref{eqmm})
we can express $\varrho(g_3)$ as
\begin{eqnarray}
\nonumber\varrho(g_3)&=&\vrh(g_1)^{-1}\vrh(g_2)\vrh(g_1)\\
\nonumber&=&\Id+(\eta-\langle e_2,\eta\rangle
e_1)(\xi+\langle\xi,e_1\rangle e_2)^T\\
&=&\Id+(\eta-\eta_2e_1)(\xi+\xi_1e_2)^T.\label{rg3}
\end{eqnarray}

\nid If we apply $\vrh$ to (\ref{EW2}) and use (\ref{eqm}) and
(\ref{rg3}), we obtain the equation
\begin{eqnarray}
\vrh(g_2)&=&\nonumber\varrho(g_3)^{-1}\vrh(g_1)\vrh(g_3)\\
&=&\nonumber\Id+\big(e_1-\langle\xi+\xi_1e_2,e_1\rangle
    \left(\eta-\eta_2e_1\right)\big)
    \big(e_2+\langle e_2,\eta-\eta_2e_1\rangle
    \left(\xi+\xi_1e_2\right)\big)^T\\
&=&\Id+\big(e_1-\xi_1(\eta-\eta_2e_1)\big)\big(e_2+\eta_2(\xi+\xi_1e_2)
    \big)^T \label{dihsss}
\end{eqnarray}

\nid Substituting this into (\ref{eqmm}) yields the equation
\begin{equation}\label{eqpmm}
(\pm\eta,\pm\xi)=\big(e_1-\xi_1(\eta-\eta_2e_1),e_2+\eta_2(\xi+\xi_1e_2)\big)
\end{equation}

\nid which is a necessary and sufficient condition for 
(\ref{eqm}) - (\ref{eqmm}) to define an $M$-homomorphism.
Now observe that the right hand side of (\ref{dihsss}) remains
unchanged if we replace $(\xi,\eta)$ by $(-\xi,-\eta)$, hence if we
choose the sign appropriately (\ref{eqpmm}) is equivalent to the
following two equations
\begin{eqnarray}
\label{EW4}(1+\xi_1)\eta&=&(1+\xi_1\eta_2)e_1\\
\label{EW5}(1-\eta_2)\xi&=&(1+\xi_1\eta_2)e_2.
\end{eqnarray}

\nid We first want to find all solutions to 
these equations with $\<\eta,\xi\>=0$ 
and $\eta,\xi$ primitive.
Let us look at (\ref{EW4}) to begin with. We distinguish the two cases
whether the sides of the equation are or are  not zero. If they are
not we can argue as follows: $\eta$ is a multiple of $e_1$,
hence being primitive it must equal $\eta=\pm e_1$. 
But then $\eta_2=0$ and (\ref{EW5}) says $\xi= e_2$. Therefore,
$\xi_1=0$ and going back to
(\ref{EW4}) this determines the sign of $\eta=e_1$. So 
$\varrho(g_2)=\Id+e_1e_2^T=A$ which gives us the 
trivial $M$-re\-pre\-sen\-ta\-tion $\vrh_{\mathrm{trivial}}$ 
(trivial in the sense
that it factors through 
$G(\grk)^{\mathrm{ab}}$). So now we assume that both sides in (\ref{EW4})
are zero. This means $\xi_1=-1$ and $\eta_2=1$. 
Then because both $\eta$ and $\xi$ have a coefficient $\pm1$ 
the primitivity condition places no further restriction on what the
other $\eta_j,\xi_j\in\ZZ$ can be, and the only other equation to
be satisfied is
\[0=\langle\eta,\xi\rangle=-\eta_1+\xi_2+\eta_3\xi_3.\]

\nid So let $a,b,c\in\ZZ$ be any integers and let us make an Ansatz
$\eta_1=a$, $\eta_3=b$, $\xi_3=c$. 
Having fixed those coefficients, $\xi_2=a- bc$ is uniquely
determined. This yields
\begin{equation}
\label{g2}
\varrho(g_2)=\Id+\left(\begin{array}{c}a \\ 1\\b \end{array}\right)
\Big(- 1,a- bc,c\Big)=\left(\begin{array}{ccc}
1- a &a^2- abc&ac \\
-1&a-bc+1&c\\
- b &ab- b^2c &bc+1\end{array}\right).
\end{equation}

\nid Finally it is easy to determine the image of $g_3$ from
(\ref{EW1}) or (\ref{rg3}); it is given by
\begin{equation}
\vrh(g_3)=\label{g3}\left(\begin{array}{ccc}
2-a &1-2a+a^2-abc+bc&ac-c \\
-1&a-bc&c\\
- b &ab- b-b^2c &bc+1\end{array}\right).
\end{equation}

Now conversely, using the Wirtinger relations
(\ref{EW1}), (\ref{EW2}) 
and (\ref{EW3}) one can easily check that for all $a,b,c\in\ZZ$
(\ref{g2}), (\ref{g3}) extends to a homomorphism 
$\varrho_{abc}\colon G(\grk)\to\SL(3,\ZZ)$ that sends all meridians
into the conjugacy class of $A$. \\

Are any of these $M$-homomorphisms equivalent? We check that the matrix
\[C=\left(\begin{array}{ccc}1&a&c\\0&1&0\\0&b&1\end{array}\right) \]

\nid commutes with $A$, and 
$C^{-1}\varrho_{abc}C=\varrho_{000}$, and all the $M$-homomorphisms
$\varrho_{abc}$ are equivalent via conjugation to the $M$-homomorphism
$\varrho_{000}$ sending
\begin{equation}\label{evrh}
g_1\longmapsto\left(\begin{array}{ccc}1&1&0\\0&1&0\\0&0&1\end{array}
\right),\;\;\; 
g_2\longmapsto\left(\begin{array}{ccc}1&0&0\\-1&1&0\\0&0&1\end{array}
\right),\;\;\;  
g_3\longmapsto\left(\begin{array}{ccc}2&1&0\\-1&0&0\\0&0&1\end{array}
\right).
\end{equation}

\nid On the other side, $\varrho_{000}$ is certainly not equivalent to
the trivial (cyclic) $M$-ho\-mo\-mor\-phism $\varrho_{\mathrm{trivial}}$
with $\varrho_{\mathrm{trivial}}(g_j)=A$ for $j=1,2,3$ since
$\varrho_{\mathrm{trivial}}(G(\grk))\cong\ZZ$ while
the matrices $A=\varrho_{000}(g_1)$ and $A'=\varrho_{000}(g_2)$ 
are well-known to generate a subgroup of $\SL(3,\ZZ)$ isomorphic to
$\SL(2,\ZZ)$.
\epf

This formalism can be used to determine the $M$-representations of more 
complicated knots (see also \refas{A2}).

% ^ 

\sec{S4}{Topology of $T^3$-fibrations} 

\ssec{SS4}{Some general facts about knots and knot groups}

The following facts are well-known and only listed for reference. See for 
instance \cite{KB} and \cite{BZ} for proofs and more details.

\prop{L976}{Let $\grk\subset S^3$ be a link with $\mu$ components and
   $B_0=S^3\setminus\grk$.
Then the integral (co-) homology of $B_0$ is given by
\begin{center}
\begin{tabular}{|c||c|c|c|c|}
\hline
$p$&$0$&$1$&$2$&$3$\\
\hline\hline
$\rmH^p(B_0,\ZZ)$&$\ZZ$&$\ZZ^\mu$&$\ZZ^{\mu-1}$&$0$\\
$\rmH_p(B_0,\ZZ)$&$\ZZ$&$\ZZ^\mu$&$\ZZ^{\mu-1}$&$0$\\
\hline
\end{tabular}
\end{center}
}\eprop

In order to compute the (co-) homology of a $T^3$-bundle $f_0\colon 
X_0\to S^3\setminus
\grk$ we will often use the Leray spectral sequence of a fibre bundle, 
whose $E_2$-term (in the case of cohomology) is a double 
complex $(E_2^{p,q})_{p,q}$ of modules 
\begin{equation}\label{rew}
E_2^{p,q}=\rmH^p(S^3\setminus \grk,
R^qf_{0*}\ZZ).\end{equation}

\nid These modules can, under certain circumstances, be identified with 
the group cohomology modules of $G(\grk)$ twisted by $\vrh$ and can thus 
be determined algebraically from the abstract representation $\vrh$. We
explain this more precisely.

\defe{Dirred}{We call a link $\grk\subset S^3$ \underline{irreducible} 
if every embedding $\iot\colon S^2\hookrightarrow S^3\setminus \grk$ is the 
restriction to $\d\oB^3=S^2$ of an embedding $\iot'\colon\oB^3
\hookrightarrow
S^3\setminus \grk$ of a closed three-ball.}\edefe

\nid With this terminology, we can quote the following well-known theorem.
\th{TP}{{\rm{(Papakyriakopoulos \cite{Papa})}} If $\grk\subset S^3$ is
  an irreducible link, then $B_0=S^3\setminus\grk$
  is an Eilenberg-MacLane space $K(\pi_1(B_0),1)=K(G(\grk),1)$.}\eth

Now $H^q(T^n,\ZZ)=\wedge^qH^1(T^n,\ZZ)$ and therefore $R^qf_*\ZZ$ is the 
local system with monodromy $\wedge^q\vrh$.
So if $\grk$ is irreducible,
the $E_2$-term of the Leray spectral sequence (\ref{rew})
can be 
interpreted as the double complex whose $(p,q)$-term is the 
group cohomology module
\begin{equation}\label{rewwwwww}
E_2^{p,q}\cong \rmH^p(G(\grk),\wedge^q\vrh).
\end{equation}

We can use the powerful machinery of group cohomology to compute these 
modules. It will be particularly helpful that we do not have to mess about 
with any explicit representing cochains. 
The case of homology is similar.
If $G$ operates on a $\ZZ$-module $\Lam$ via $\vrh$ we also
write $H_*(G,\Lam)$ for $H_*(G,\vrh)$.
If we use the ``bar resolution'' to define the (co-) homology of a group $G$ 
with coefficients in $\Lam$, the following formulas for the 
low-dimensional (co-) homology modules are immediate.

\prop{L54}{Let $G$ be a group and $\Lam$ a $G$-module. Then 
\begin{eqnarray*}
\rmH_0(G,\Lam)&=&\frac{\Lam}{\langle g\xi-\xi\;|\;g\in G, \xi\in
  \Lam\rangle}\;=\;\Lam_G\\ 
\rmH^0(G,\Lam)&=&\{\xi\in \Lam\;|\; \forall g\in G\colon g\cdot\xi=\xi
  \}\;=\;\Lam^G\\ 
\rmH^1(G,\Lam)&=&\frac{\mathrm{Der}(G,\Lam)}{\mathrm{Der}_0(G,\Lam)}
\end{eqnarray*}}\eprop

Here $\Lam^G$ is the space of {\em $G$-invariants}, the largest 
submodule of $\Lam$ 
on which $G$ acts trivially, while $\Lam_G$ is the space of 
{\em $G$-coinvariants}, 
the largest quotient of $\Lam$ on which $G$ acts trivially.

Computing the higher degree 
homology or cohomology of a group can be rather difficult, 
even for everyday groups. However, the computation of $H_*(G,\Lam)$ and
$H^*(G,\Lam)$ is quite easy if $G$ is infinite cyclic. 

\prop{Pcyc}{Suppose $G$ is a free cyclic group and $\Lam$ is a $G$-module. 
Then 
\[\begin{array}{lclcl}
\rmH^0(G,\Lam)&\cong&\rmH_1(G,\Lam)&\cong&\Lam^G\\
\rmH^1(G,\Lam)&\cong&\rmH_0(G,\Lam)&\cong&\Lam_G\\
\rmH^i(G,\Lam)&\cong&\rmH_i(G,\Lam)&\cong&\{0\}\;\;\;\;\;\mathrm{for\;\;}i>1.
\end{array}\]
}\eprop

In case a group $G$ is an amalgamated sum of groups with known (co-) homology 
the following
Mayer-Vietoris sequence can be used to compute the (co-) homology of $G$ 
from that of its factors.

\th{TEX}{{\rm{(Lyndon \cite{Ly}, Swan \cite{Sw}) }} Let $G=G_1*_AG_2$
  be a free product with amalgamated subgroup 
  $A$, and let $\vrh\colon G\to\Aut_\ZZ(\Lam)$ be a representation of $G$
  on a $\ZZ$-module $\Lam$, making $\Lam$ into a $\ZZ G$-module. 
  Then the following naturally defined Mayer-Vietoris sequences are exact:  
\[\begin{array}{c}
\cdots\to\rmH_{n+1}(G,\Lam)\to\rmH_n(A,\Lam)\to\rmH_n(G_1,\Lam)
\oplus\rmH_n(G_2,\Lam)\to 
\rmH_n(G,\Lam)\to\cdots\\[0.8ex]
\cdots\to\rmH^{n-1}(A,\Lam)\to \rmH^n(G,\Lam)\to\rmH^n(G_1,
\Lam)\oplus\rmH^n(G_2,\Lam)\to 
\rmH^n(A,\Lam)\to\cdots
\end{array}\]}\eth

 The fact that if $G(\grk)$ is the group of a knot or irreducible link
 $\grk\subset S^3$ then $S^3\setminus\grk$ is a $K(G(\grk),1)$ implies 
 a strong interplay between the topology of a $T^3$-fibration over $S^3$ 
 whose discriminant locus is a knot or an irreducible link
 and the representation theory of the corresponding knot group.

\prop{C1}{The group of a knot or irreducible link has no cohomology with
  coefficients in degree $>2$ for any re\-pre\-sen\-ta\-tion.}\eprop

\pf
If $G(\grk)$ is the group of a knot 
or irreducible link $\grk\subset S^3$ then $S^3\setminus\grk$ is a 
$K(G(\grk),1)$ space by \reft{TP}, so the group (co-) homology of $G(\grk)$
is identified with the singular (co-) homology of $S^3\setminus \grk$ which is 
a smooth open three-manifold of cohomological dimension $\cd(S^3\setminus
\grk)\leq2$.\epf

\ssec{SS42}{Fundamental group and characteristic classes}

We now show how to compute the topological invariants of an affine
$T^3$-fibration of the type considered in this paper
from its monodromy representation. 

\th{Psc}{Suppose $\grk$ is an irreducible link and $b\in S^3\setminus \grk$.
  Let $G(\grk)=\pi_1(S^3\setminus\grk,b)$ and $\vrh\colon G(\grk)\to 
\SL(3,\ZZ)$
  an $M$-re\-pre\-sen\-ta\-tion. If
\begin{equation}\label{eqvan}
  \rmH_0(G,\vrh)=
  \frac{\rmH_1(F_b,\ZZ)}{\{\varrho(g)\xi-\xi\;|\;\xi\in\rmH_1(F_b,\ZZ),\;g\in 
  G(\grk)\}} =\{0\},
\end{equation}

  the total space $X$ of the good $T^3$-fibration (with section) $f\colon
  X\to S^3$ associated to $\vrh$
  is simply connected.}\eth

\pf Let $b\in S^3\setminus \grk$ 
and $\iot\colon F_b\to X$ the inclusion of the smooth
fibre over $b$ in $X$.
We claim that (\ref{eqvan}) implies that for
every loop $\gam$ in $F_b$ with base point $\sig(b)$, the push-forward
$\iot_*\gam$ is null-homotopic in $X$. The reason is, geometrically speaking,
that ``$H_1(F_b,\ZZ)$ is generated by vanishing cycles''. We will explain 
this more precisely now. Note that $G(\grk)$ is generated by
finitely many meridians $g_1\ddd g_l$ and 
\[\vrh(g_ig_j)\xi-\xi=[\vrh(g_i)(\vrh(g_j)\xi)-
\vrh(g_j)\xi)]+[\vrh(g_j)\xi-\xi],\]

\nid thus by induction on the length of words in $g_1\ddd g_l$ we
conclude from (\ref{eqvan}) that every $\xi\in H_1(F_b,\ZZ)$ can be given 
in terms of the $g_i$ by an expression
\begin{equation}\label{eqvan2}
\xi=\sum\limits_{i=1}^r[\varrho(g_{k_i})\xi_i-\xi_i],\;\;\;\;k_i\in\{1\ddd
l\},\;\;\xi_i\in\rmH_1(F_b,\ZZ).\end{equation}

\nid The image of the map 
\[(\varrho(g_i)-\Id)\colon
\rmH_1(F_b,\ZZ)\to \rmH_1(F_b,\ZZ)\]

\nid is a $\ZZ$-module $E_i\subset
H_1(F_b,\ZZ)$ with $\rk(E_i)=1$, and by (\ref{eqvan}) and (\ref{eqvan2}),
\[E_1+\dots+ E_l=\rmH_1(F_b,\ZZ).\]

\nid In our local model (\ref{Diag}), there is a fibre-preserving 
$S^1$-action in a neighbourhood of every singular fibre which
is free outside the singular points, and whose fixed points are
precisely the critical points. If $x\in F_b$, the orbit of $x$ under
the $S^1$-action represents a 
non-trivial homology class $[\eps_i]\in 
H_1(F_b,\ZZ)$ which generates
$E_i$ ($[\eps_i]$ is the ``vanishing cycle'' for $g_i$). 
Now let $[\gam]\in\pi_1(F_b,\sig(b))$, then since
$H_1(F_b,\ZZ)\cong\pi_1(F_b,\sig(b))$ we can write the representing loop
\[\gam=\sum\limits_{i=1}^l r_i\eps_i,\;\;\;\;\;r_i\in\ZZ,\] 

\nid thus realizing $[\gam]$ as a ``linear combination'' of vanishing
cycles. But now the
collapsing $S^1$-action yields a homotopy from $\gam$ to a constant
loop, and 
\[\iot_*(\pi_1(F_b,\sig(b)))=\{1\}\subset\pi_1(X,\sig(b)).\]
Now let $\gam\colon S^1\to X$ be an arbitrary (say, piecewise linear or
smooth) loop. Since $f^{-1}(\grk)$ is contained in a finite union of
codimension $2$ submanifolds of $X$, we can move $\gam$ slightly
inside its homotopy class to make it disjoint from $C=f^{-1}(\grk)$. It
then represents an element of $\pi_1(X\setminus C,\sig(b))$.
Now we have the following homotopy exact sequence (omitting basepoints):
\[\begin{CD}\pi_2(S^3\setminus
  \grk)@>>>\pi_1(F_b)@>>>\pi_1(X\setminus C)
  @>>>\pi_1(S^3\setminus
  \grk)@>>>\pi_0(F_b)\\ 
@V{\cong}VV @V{\cong}VV @V{\cong}VV @V{\cong}VV
  @V{\cong}VV \\
\{0\}@>>>\rmH_1(F_b,\ZZ)@>>>\rmH_1(F_b,\ZZ)\rtimes G(\grk) @>>> G(\grk) @>>> 
\{*\}
\end{CD}\]

\nid Because $X\to S^3$ has a section $\sig$ this sequence splits and
we can write  
\[\gam=(\iot_*\gam_1)(\sig_*\gam_2)\]

\nid with $[\gam_1]\in\pi_1(F_b,\sig(b))$ and $[\gam_2]\in\pi_1(S^3\setminus
\grk,b)$. However, we just said that under the inclusion
$\pi_1(X\setminus C,\sig(b))\to\pi_1(X,\sig(b))$ the element
$[\iot_*\gam_1]\mapsto1\in\pi_1(X,\sig(b))$ and 
$[\sig_*\gam_2]\mapsto1\in\pi_1(X,\sig(b))$ because $S^3$ is simply
connected. Thus 
\[\pi_1(X,\sig(b))=\{1\}\]
and $X$ is simply connected. This proof works more generally 
for $T^n$-fibrations
($n\geq2$) with singularities of the form (\ref{Diag}) as long as 
(\ref{eqvan}) holds. 
\epf

\prop{Pori}{Every good $T^3$-fibration $X\to S^3$ degenerating over
  a link $\grk$ with monodromy 
  representation $\vrh\colon G(\grk)\to\SL(3,\ZZ)$  is orientable.}\eprop 

\pf Let $X_0=X\setminus f^{-1}(\grk)$ be the smooth part of the fibration 
and $B_0=S^3\setminus\grk$. Then 
\[TX_0=\calV\oplus f^*TB_0\]

\nid where $\calV$ is the vertical bundle of the $T^3$-bundle 
$f_0\colon X_0\to B_0$. Now 
\[\calV=B_0\times_\vrh\RR^n,\]
hence we get
\[\Lam^{6}TX_0=(B_0\times_{\wedge^3\vrh}\RR)\otimes f^*(\Lam^3TS^3)|_{B_0}.\]

\nid This is trivial because $\vrh(G(\grk))\subset\SL(3,\ZZ)$ and so 
$\wedge^3\vrh$ is trivial, as is $TS^3$, hence $X_0$ is orientable. 
But $\dim(X\setminus X_0)=4=\dim(X)-2$
and so $X$ itself is also orientable.\epf

\th{Cff}{Let $f\colon X\to S^3$ be a good
  $T^3$-fibration with section $\sig$ and $X^\#=X\setminus \Crit(f)$. Then
  the vertical bundle $\calV\subset TX^\#$ and $TX^\#$ itself 
  are trivial.}\eth

\pf The linear structure on the $T^3$-bundle $X_0\to S^3\setminus
\grk$ extends to the smooth part of the singular fibres, hence there is
a short exact sequence of sheaves
\[ 0\lra R^{n-1}(f^\#)_*\ZZ\lra \sig^*\calV\lra X^\#\to0\]
\nid and hence 
\[\calV\cong f^*(\sig^*\calV)\]
\nid s trivial as a pull-back of a bundle on the base $S^3$. (Every
vector bundle on $S^3$ is trivial since $\pi_2(\GL(n,\RR))=\{0\}$). Moreover,
\[TX^\#=\calV\oplus f^*TS^3|_{X^\#}\]
and $TS^3$ is trivial, hence $TX^\#$ is trivial, too. \epf

\th{TSW}{The underlying six-manifold $X$ of a good $T^3$-fibration
  $f\colon X\to S^3$ has vanishing second
  Stiefel-Whitney class, $w_2(X)=0$.}\eth

\pf (see \cite{G3})
By \refp{Pori}, a good $T^3$-fibration is orientable. 
Since $\Crit(f)\subset X$ is of codimension $\geq 4$, we have
$w_2(X)=w_2(X^\#)$. But $TX^\#$ is trivial by
\reft{Cff},
thus $w_2(X^\#)=0$.
\epf

\th{TPC}{In a six-manifold $X$ carrying a good $T^3$-fibration, the critical
  locus represents the first Pontryagin class: 
\[c\,[\Crit(f)]=p_1(X)\in H^4(X,\ZZ)\]
  for some constant $c\in\ZZ$ under Poincar{\'e}-duality.\label{poincare}}\eth

\pf This theorem is also contained in \cite{G3}, but since we are in the
simpler situation of good $T^3$-fibrations, we can avoid the use 
of $K$-theory.

\pid Let $f\colon X\to
S^3$ be a good $T^3$-fibration. Then the bundle $TX^\#$ is
trivial. Choose a trivialisation of $TS^3$, this gives us a preferred
trivialisation of $TX^\#$ which in turn
defines a relative Pontryagin class
$p_{\mathrm{rel}}(X,X^\#)\in H^4(X,X^\#,\ZZ)$. There is a short exact
sequence
\[\rmH^3(X,\ZZ)\stackrel{\phi}{\lra}\rmH^3(X^\#,\ZZ)\stackrel{\del}{\lra}
\rmH^4(X,X^\#,\ZZ)\stackrel{\psi}{\lra }\rmH^4(X,\ZZ).\]
\nid The image $p_1(X)=\psi(p_{\mathrm{rel}}(X,X^\#))$ is the first
Pontryagin class of $X$. (It is independent
of the chosen trivialisation of $TS^3$ because any other choice 
is given by a map $t\colon S^3\to \SO(3)$ and changes
$p_{\mathrm{rel}}(X,X^\#)$ by $\del((f^\#)^*t^*\ome)$ where $\ome\in
H^3(\SO(3))$ is the canonical degree-$3$ class. But $(f^\#)^*(t^*\ome)=
\phi(f^*(t^*\ome))$ and hence $\del((f^\#)^*t^*\ome)=0$.) But now
$H^4(X,X^\#,\ZZ)\cong H_2(X\setminus X^\#,\ZZ)\cong\ZZ$ 
and $p_1(X)$ is represented
by some multiple of the generator $[\Crit(f)]\in H_2(X\setminus X^\#,\ZZ)$. 
\epf

\th{TECH}{For any good $T^3$-fibration $f\colon X\to B$ the
  Euler-char\-ac\-ter\-is\-tic vanishes, $\chi(X)=0$.}\eth

\pf For a fibre bundle $f\colon X\to B$ with fibre $F$ we
have the formula 
$\chi(X)=\chi(F)\chi(B)$. Write the total space of our
$T^3$-fibration as $X=X_0\cup X_1$ where $X_0\to B\setminus\grk$ is
the smooth part of the fibration
and $X_1\to\grK$ is the singular fibration over an open solid
torus. Then
\[\chi(X)=\chi(X_0)+\chi(X_1)-\chi(X_0\cap X_1).\]

\nid Now $\chi(X_0)=\chi(T^3)\chi(B\setminus\grk)=0$, and $X_1$ is homotopy
equivalent to a fibre bundle $X_1'\to S^1$ with fibre
$F'\cong I_1\times S^1$, 
hence $\chi(X_1)=\chi(F')\chi(S^1)=0$. Finally $\chi(X_0\cap X_1)=0$ because
the intersection retracts onto a compact smooth five-manifold (a
$T^3$-bundle over $T^2$). Thus $\chi(X)=0$. \epf

% ^ 

\sec{S5}{Examples}

We study the $T^3$-fibrations associated to the simplest non-trivial 
$M$-representations of torus knot groups
discussed in \refs{S4}, namely those $M$-representations whose image is 
a subgroup of $\SL(3,\ZZ)$ isomorphic to $\SL(2,\ZZ)$. Geometrically 
this means that we first construct a $T^2$-fibred five-manifold 
$Y\to S^3$ and then consider $X=\TT^1\times Y$ or, more generally, 
(if we do not insist on the existence of a section) 
principal-$\TT^1$-bundles $X\to Y$.

\ssec{S4SS1}{The trefoil $\grt(2,3)$ and $S^3\times S^3$}

Let $\grk=\grt(2,3)\subset S^3$ be the trefoil knot.
We now show that the non-trivial $M$-representation of the trefoil group
\[\vrh\colon G(\grk)\to\SL(2,\ZZ)\subset\SL(3,\ZZ)\]
\nid which was given in \reft{P2345610} (for an explicit formula 
see (\ref{evrh})) gives us $T^3$-fibrations on spaces such as
$S^1\times S^2\times S^3$ and 
$S^3\times S^3$.\\

\pid We begin with a short detour which offers a geometric description of the
$T^3$-fibration associated to $\vrh$.
Recall the classical fact (see \cite{Mi1}) that 
\begin{equation}\label{ios}
S^3\setminus \grk\cong \frac{\SL(2,\RR)}{\SL(2,\ZZ)}.
\end{equation}
\nid The right hand side is naturally the space of lattices 
of unit volume in $\CC$. To any lattice $\Lam\subset\CC$ we can 
associate the unique pair $(g_2(\Lam),g_3(\Lam))$ of complex numbers 
\begin{equation}\label{tyu}
g_2(\Lam)=60\sum_{\ome\in\Lam\setminus\{0\}}\frac{1}{\ome^4},\;\;\;\;\;\;\;
g_3(\Lam)=140\sum_{\ome\in\Lam\setminus\{0\}}\frac{1}{\ome^6}
\end{equation}
so that the elliptic curve
$\CC/\Lam$ is equivalent to an elliptic curve in $\CC P^2$
defined by the Weierstra{\ss} equation 
\begin{equation}\label{eqeqeqeqeq}
4x^3-g_2(\Lam)xz^2-g_3(\Lam)z^3-y^2z=0,
\end{equation} 
where $[x:y:z]$ are homogeneous coordinates on $\CC P^2$.
If we rescale the lattice \mbox{$\Lam\to t\Lam$}
by some $t>0$, its invariants transform by 
\[g_2(t\Lam)=t^{-4}g_2(\Lam),\;\;\;g_3(t\Lam)=t^{-6}g_3(\Lam)\]
as (\ref{tyu}) shows. So if $\Lam_0=\ZZ\oplus i\ZZ\subset\CC$ is 
the standard lattice (of Gau{\ss} integers)
there exists a unique function $t\colon \SL(2,\RR)\to\RR^+$ so that  
$\Lam_M=t(M)M\Lam_0$ satisfies
\[|g_2(\Lam_M)|^2+|g_3(\Lam_M)|^2=1.\]
for all $M\in\SL(2,\RR)$,
where the action of $\SL(2,\RR)$ on the complex 
plane comes from the standard isomorphism $\RR^2\cong\CC$.
Obviously $g_i(\Lam_M)=g_i(\Lam)$ for $i=2,3$ if and only if $M\in\SL(2,\ZZ)$.
Given $M\in\SL(2,\RR)$, the isomorphism (\ref{ios}) sends the 
$\SL(2,\ZZ)$ orbit of $M$ to $(g_2(\Lam_M),g_3(\Lam_M))\in S^3
\subset\CC^2$.
A pair $(g_2,g_3)\in S^3$ can arise 
in this manner
if and only if the discriminant satisfies
\begin{equation}\label{disc}
\Del=27g_3^2-g_2^3\neq0.
\end{equation}
It is well-known \cite{Mi2} that (\ref{disc}) describes the 
complement of a trefoil in $S^3$.\\

\pid Using this isomorphism 
let us now define a ``tautological'' $\ZZ^2$-bundle $\calL\to S^3\setminus
\grk$ whose fibre at $\Lam=(g_2,g_3)\in S^3\setminus \grk$ is 
$\Lam$ itself. In other words, if we embed
\begin{eqnarray*}
\vph\colon \SL(2,\RR)\times\Lam_0&\hookrightarrow&\SL(2,\RR)\times\CC\\
(M,\xi)&\mapsto& (M,t(M)M\xi),
\end{eqnarray*}
then by the isomorphism (\ref{ios}) $\calL$ is the quotient
\[\calL=\frac{\IM(\vph)}{\SL(2,\ZZ)}\]

where $\SL(2,\ZZ)$ acts on $\SL(2,\RR)$ by right-multiplication 
and on $\CC$ by 
its real standard representation.\\

\pid Now from (\ref{evrh}) there is a short exact sequence
\[\{1\}\longrightarrow 2\calZ(G)\longrightarrow
G\stackrel{\vrh}{\longrightarrow}\SL(2,\ZZ)\longrightarrow \{1\}.\]
The centre $\calZ(G)\cong\ZZ$ is infinite cyclic and $2\calZ(G)=\ker(\vrh)$ 
is the subgroup of index two in $\calZ(G)$. Thus the $M$-representation 
$\vrh$ descends
to a representation $\bar{\vrh}$ of the quotient 
$G(\grk)/2\calZ(G)\cong\SL(2,\ZZ)$, and it is easy to check explicitly that 
$\bar{\vrh}$ is the standard representation of $\SL(2,\ZZ)$.
Let $\widehat{\SL}(2,\RR)\to\SL(2,\RR)$ be the universal cover. Then
\begin{eqnarray*}
Y_0&=&\widehat{S^3\setminus\grk}\times_\vrh\TT^2\\
&=&\widehat{\SL}(2,\RR)\times_\vrh\TT^2\\
&=&\SL(2,\RR)\times_{\bar{\vrh}}\TT^2\\
&=&\frac{\SL(2,\RR)\times_{\bar{\vrh}}\CC}{\calL}.
\end{eqnarray*}
This is a $T^2$-bundle
whose fibre at $\Lam=(g_2,g_3)\in S^3\setminus\grk$ is
$\CC/\Lam$. Hence we can see that there is a natural embedding 
$Y_0\subset \CC P^2\times (S^3\setminus\grk)$ given by the equation
 \begin{equation}\label{1eqeqeqeqeq}
4x^3-g_2xz^2-g_3z^3-y^2z=0\end{equation}
with $[x:y:z]\in\CC P^2$ and $(g_2,g_3)\in S^3$. Note that
the $T^2$-bundle $f\colon Y_0\to S^3\setminus\grk$ also has a 
section $\sig_0\colon S^3\to X$ mapping
\[(g_2,g_3)\in S^3\longmapsto ([0:1:0], (g_2,g_3))\in X,\]
where $\sig(g_2,g_3)$ is the point at infinity of $f^{-1}(g_2,g_3)$.
By construction its
  monodromy is given
  by the unique non-trivial $M$-re\-pre\-sen\-ta\-tion $\varrho\colon
  G(\grt{(2,3)})\to \SL(2,\ZZ)$.
However, equation (\ref{1eqeqeqeqeq}) 
is well-defined on $\CC P^2\times S^3$ and thus 
defines a natural compactification $Y_0\subset Y\subset \CC P^2\times S^3$.
The section extends to a section $\sig\colon S^3\to Y$.

\prop{Psmo}{$Y$ is a smooth compact $T^2$-fibred five-manifold.}\eprop

\pf Use inhomogeneous coordinates to write $Y$ as the vanishing locus
of a smooth function $f$ with $0$ as a regular value. 
The unitary condition $g_2\bar{g}_2+g_3\bar{g}_3=1$ ensures that
there are no critical points in $f^{-1}(0)$.
The implicit function theorem then shows that $Y$ is smooth, and it is
obviously compact.\epf

From \reft{TSW} we know that $Y$ is spin, but we can also see this directly
in this situation.

\th{TSpin}{$Y$ is a spin manifold, that is the second Stiefel-Whitney
  class vanishes, \mbox{$w_2(Y)=0$}.}\eth

\pf Let $\iot\colon Y\hookrightarrow\CC P^2\times S^3$ be the inclusion and let
\[\varpi\in\rmH^2(\CC P^2
\times S^3,\ZZ/2)\cong\rmH^2(\CC P^2,\ZZ/2)\cong\ZZ/2\]
be the non-trivial degree 2 cohomology class.
Let $\calN\subset T(\CC P^2\times S^3)|_Y$ be the normal bundle of
the submanifold $Y$, then
\[w_2(Y)+w_2(\calN)=\iot^*w_2(\CC P^2\times S^3) =\iot^*(\mathrm{pr}_1^*w_2
(\CC P^2)) = \iot^*\varpi.\]
But since the equation (\ref{eqeqeqeqeq}) is of odd degree $3$ we get
$\calN=\iot^*p_1^*\calO(-3)$ and 
\[w_2(\calN)=-3\iot^*\varpi=\iot^*\varpi,\]
\nid and thus $w_2(Y)=0$ and $Y$ is spin.\epf

We would now like to identify this manifold. We first compute its
integral homology. We define:
\begin{equation}\label{notat}
\begin{array}{lclclcl}
B_0&=&S^3\setminus \grk{}{},&&Y_0=f^{-1}(B_0),\\
B_1&=&\calN(\grk{}{}),&&Y_1=f^{-1}(B_1).
\end{array}.\end{equation}
\nid As usual, $\calN(\grk{}{})$ denotes a tubular neighbourhood
of the knot.

\th{pcht23}{The integral homology of the five-manifolds
$Y_0$, $Y_1$ and $Y_0\cap Y_1$ is as follows:
\begin{center}
\begin{tabular}{|c||c|c|c|}
\hline
$p$&$\rmH_p(Y_0)$&$\rmH_p(Y_1)$&$\rmH_p(Y_0\cap Y_1)$\\
\hline\hline
$0$ &$\ZZ$            &$\ZZ$            &$\ZZ$\\
$1$ &$\ZZ$            &$\ZZ\oplus\ZZ/2$ &$\ZZ^2\oplus\ZZ/2$\\
$2$ &$\ZZ\oplus\ZZ/2$ &$\ZZ$            &$\ZZ^2\oplus\ZZ/2$\\
$3$ &$\ZZ$            &$\ZZ$            &$\ZZ^2$\\
$4$ &$0$              &$0$              &$\ZZ$\\
$5$ &$0$              &$0$              &$0$\\
\hline
\end{tabular}
\end{center}}\eth

\nid\pf {\it Homology of $Y_0$}. We can compute $H_*(Y_0,\ZZ)$ from the
Leray spectral sequence of the fibre bundle whose $E^2$-term we now show 
is the following:
\begin{equation}\label{eqlelel}
\begin{array}{cccl}
\;\ZZ\;&\;\ZZ\;&\;0\;&\;0\;\\
\;0\;&\;\ZZ/2\;&\;0\;&\;0\;\\
\;\ZZ\;&\;\ZZ\;&\;0\;&\;0\;
\end{array}\end{equation}
\nid Here $E^2_{p,q}\cong H_p(G,\wedge^q\vrh)$ by (\ref{rewwwwww}) for
$0\leq p\leq3$ and $0\leq q\leq2$. Since
$\wedge^0\vrh=\wedge^2\vrh$ is trivial we have $E^2_{p,0}\cong
E^2_{p,2}\cong H_p(B_0,\ZZ)$ and so the top and bottom rows are given
by \refp{L976}. For the middle row we have from \refp{L54}
\[E_{0,1}^2=\rmH_0(G,\vrh)=\frac{\Lam_0}{\{g\xi-\xi\;|
  \;\xi\in \Lam_0,\;\;g\in G\}}\cong\frac{\ZZ^2}{\{\varrho(g)\xi-\xi\;|
  \;\xi\in \ZZ^2,\;\;g\in G\}}.\]
\nid However, with respect to the standard basis $\Lam_0\cong \ZZ e_1\oplus
\ZZ e_2$ and meridians $g_1,g_2$ as before we have
$\vrh(g_2)e_1-e_1=-e_2$ and $\vrh(g_1)e_2-e_2=e_1$, thus
the right hand side is trivial and $H_0(G,\vrh)\cong E_{0,1}^2\cong\{0\}$.
Also $H_3(G,\vrh)\cong E_{3,1}^2\cong\{0\}$. To compute $E^2_{2,1}$
  and $E^2_{1,1}$ we use the
  Mayer-Vietoris sequence. We can write $G=\<x\>*_{\<x^2=y^3\>}\<y\>$ as the
  amalgamated sum of two cyclic groups $G_1=\<x\>$ and $G_2=\<y\>$
  along the cyclic subgroup $A\cong\<x^2\>\cong\<y^3\>$.
  Cyclic groups have no homology in degree $>1$ and the
  first term which has a chance of being non-zero is $H_2(G,\vrh)\cong
  E^2_{2,1}$, so the tail of this
  sequence looks like
\begin{equation}\label{eemves}
0\longrightarrow E^2_{2,1}\longrightarrow
  \rmH_1(A,\vrh)\stackrel{\psi_1}{\longrightarrow}
  \rmH_1(G_1,\vrh)\oplus\rmH_1(G_2,\vrh)\longrightarrow\cdots.
\end{equation}
\nid By \refp{Pcyc} and \refp{L54}
\[\begin{array}{ccccc}
\rmH_1(A,\vrh)&\cong&\rmH^0(A,\vrh)&\cong&\{\xi\in\ZZ^2\;|
\;\vrh(x^2)\xi=\xi\}\\
\rmH_1(G_1,\vrh)&\cong&\rmH^0(G_1,\vrh)&\cong&
\{\xi\in\ZZ^2\;|\;\vrh(x)\xi=\xi\}\\
\rmH_1(G_2,\vrh)&\cong&\rmH^0(G_2,\vrh)&\cong&
\{\xi\in\ZZ^2\;|\;\vrh(y)\xi=\xi\}
\end{array}\]
\nid We check from the matrices given in \reft{P2345610} that
\[\begin{array}{ccccc}
\varrho(x)&=&\vrh(g_3g_2g_1)&=&
\left(\begin{array}{rr}0&1\\-1&0\end{array}\right),\\
\varrho(y)&=&\vrh(g_3g_2)&=&
\left(\begin{array}{rr}1&1\\-1&0\end{array}\right),\\
\varrho(x^2)&=&\vrh(y^3)&=&
\left(\begin{array}{rr}-1&0\\0&-1\end{array}\right).
\end{array}\]
\nid None of the re\-pre\-sen\-ta\-tions induced from
$\vrh$ on
$A,G_1$ and $G_2$ has invariants, so
\[\rmH_1(G_1,\vrh)\cong\rmH_1(G_2,\vrh)\cong\rmH_1(A,\vrh)\cong\{0\}.\]

\nid Hence $E^2_{2,1}\cong\{0\}$. But also
$\coker(\psi_1)=\{0\}$ in (\ref{eemves})
and so $E_{1,1}^2\cong\rmH_1(G,\vrh)$ is the kernel
of the map $\psi_0$ in the short exact sequence
\begin{equation}\label{eeq23}
0\to E^2_{1,1}\to \rmH_0(A,\vrh)\stackrel{\psi_0}{\to}
  \rmH_0(G_1,\vrh)\oplus\rmH_0(G_2,\vrh)
  \to\rmH_0(G,\vrh)\cong\{0\}. \end{equation}
\nid We compute these $0$-th homology groups from \refp{L54}. Note that
\begin{eqnarray*}
B=\vrh(x)-\Id&=&
\left(\begin{array}{rr}-1&1\\-1&-1\end{array}\right),\;\;\\
B'=\vrh(y)-\Id&=&
\left(\begin{array}{rr}0&1\\-1&-1\end{array}\right),\;\;\\
B''=\vrh(x^2)-\Id=\vrh(y^3)-\Id&=&
\left(\begin{array}{rr}-2&0\\0&-2\end{array}\right).
\end{eqnarray*}
\nid Since $\det(B)=2$ and $\det(B')=1$ we deduce that
\begin{eqnarray*}
\rmH_0(G_1,\vrh)&\cong&\frac{\ZZ^2}{B\ZZ^2}\cong\ZZ/2,\\
\rmH_0(G_2,\vrh)&\cong&\frac{\ZZ^2}{B'\ZZ^2}\cong\{0\},\\
\rmH_0(A,\vrh)&\cong&\frac{\ZZ^2}{B''\ZZ^2}\cong\frac{\ZZ^2}{2\ZZ^2}\cong
\ZZ/2\times\ZZ/2.\end{eqnarray*}

\nid Thus (\ref{eeq23}) becomes a short exact sequence of groups
\[0\longrightarrow
E_{1,1}^2\longrightarrow\ZZ/2\times\ZZ/2
\longrightarrow\ZZ/2\longrightarrow0.\]

\nid The only group fitting into this sequence is
$E_{1,1}^2\cong\ZZ/2$. It is obvious that the boundary operator
$\d_2\colon E^2_{p,q}\to E^2_{p-2,q+1}$
of the $E^2$-term must be zero since all non-trivial homology is
concentrated in degree $0$ and $1$. Hence the homology $Y_0$ is the
total homology of $E^2_{p,q}$ which is the homology given in the
proposition. \\

\nid {\it Homology of $Y_1$.} The singular fibration $Y_1\to B_1$ is
homotopy equivalent to a fibre bundle
of singular elliptic curves of type $I_1$ over the circle. Note that a
longitude $\ell$ acts on $\Lam_0\cong\ZZ^2$ via
\begin{equation}\label{eqlll}
\vrh(\ell)=\left(\begin{array}{rr}-1&6\\0&-1\end{array}\right).
\end{equation}

\nid Let $F_0$ be a singular fibre. Then $H_p(F_0,\ZZ)\cong\ZZ$ for
$p=0,1,2$ and $H_p(F_0,\ZZ)=\{0\}$ for $i>2$, and (\ref{eqlll}) implies
that the monodromy around $\ell$ acts as $(-1)$ on
$H_1(F_0,\ZZ)\cong\ZZ$. Cover $S^1$ with two
open intervals $U',V'$ and let $U=Y_1|_{U'}$ and $V=Y_1|_{V'}$ so that
\[\rmH_p(U,\ZZ)\oplus \rmH_p(V,\ZZ)\cong
\rmH_p(U\cap V,\ZZ)\cong \rmH_p(F_0,\ZZ)\oplus \rmH_p(F_0,\ZZ)\cong \ZZ^2\]
\nid for $p=0,1,2$. Choosing appropriate trivialisations the map
$\psi_p$ in the
Mayer-Vietoris sequence
\begin{equation}\label{eeqmvv}
\cdots \to \rmH_{p+1}(Y_1,\ZZ)\to \rmH_p(U\cap V,\ZZ)
\stackrel{\psi_p}{\to}\rmH_p(U,\ZZ)\oplus \rmH_p(V,\ZZ) \to
\rmH_{p}(Y_1,\ZZ)\to \cdots \end{equation}
\nid takes the form $(a,b)\in\ZZ^2\mapsto (a-b, a-\mu_pb)\in\ZZ^2$
where $\mu_p$ is the monodromy action on $H_p(F_0,\ZZ)$. We have
$\mu_0=\mu_2=1$ and $\mu_1=-1$.
Because $H_*(F_0,\ZZ)$ is torsion free and all modules are Abelian,
the Mayer-Vietoris sequence
splits at all of its boundary operators, and
\begin{equation}\label{eqkck}
\rmH_p(Y_1,\ZZ)=\ker(\psi_{p-1})\oplus\coker(\psi_{p})\end{equation}
for all $p$. Let $\psi^\pm\colon\ZZ^2\to\ZZ^2$ be the map
$(a,b)\mapsto(a-b,a\mp b)$, then
\begin{eqnarray*}
\ker(\psi^+)&=&\diag(\ZZ^2)\cong\ZZ,\\
\coker(\psi^+)&=&\frac{\ZZ^2}{\diag(\ZZ^2)}\cong\ZZ,\\
\ker(\psi^-)&=&\{(a,b)\in\ZZ^2\;|\;a+b=a-b=0\}\cong\{0\},
\;\;\;\;\;{\mathrm{(no\;}}2 {\mathrm{\;torsion)}}\\
\coker(\psi^-)&=&\frac{\ZZ^2}{\{(a-b,a+b)\;|\;a,b\in\ZZ\}}\cong\ZZ/2.
\end{eqnarray*}
\nid Now $\psi_0=\psi_2=\psi^+$, $\psi_1=\psi^-$ and $\psi_p=0$ for
$p>2$. Together with
(\ref{eqkck}) this implies the claim.  \\

\nid {\it Homology of $Y_0\cap Y_1$.} We compute again the Leray spectral
sequence for the $T^2$-bundle over $B_0\cap B_1\simeq T^2$. We claim
that the $E^2$-term has the following form:
\begin{equation}\label{eqlelelel}
\begin{array}{ccl}
\;\ZZ\;  &\;\ZZ^2\;&\;\ZZ\;\\
\;\ZZ/2\;&\;\ZZ/2\;&\;0\;\\
\;\ZZ\;  &\;\ZZ^2\;&\;\ZZ\;
\end{array}\end{equation}
\nid Again, the top and bottom rows are just the integral homology of
$B_0\cap B_1\simeq T^2$. For the middle row
note that the inclusion $B_0\cap B_1\hookrightarrow B_0$ is an
inclusion of $K(\pi_1,1)$-complexes corresponding to the embedding
$G'=\<m,\ell\>\to G$ of the subgroup of $G$ generated by a longitude
$\ell$ and a meridian $m$ into $G$. Hence let $\vrh'=\vrh|G'$. Then
$E^2_{p,1}\cong H_p(G',\vrh')$ for $p=0,1,2$. As before we compute
$H_0(G',\vrh')$ from \refp{L54} and equations (\ref{eeaa}) and
(\ref{eqlll}) which give explicitly the actions of $m$ and $\ell$ on
$\Lam_0$ with respect to the standard trivialisation of $\Lam_0\cong\ZZ^2$. We
find
\[E^2_{0,1}\cong\rmH_0(G',\vrh')
\cong\ZZ/2.\]
\nid Note that $\vrh\sim\vrh^*$ because the representation takes values in 
$\SL(2,\ZZ)$ and
\[E^2_{2,1}\cong\rmH_2(G',\vrh')
  \cong\rmH^0(G',(\vrh')^*)\cong
  \rmH^0(G',\vrh')\cong\{0\}\]
\nid by Poincar{\'e}-duality ((\ref{eqlll}) shows that
$\vrh'$ has no invariants).
Finally in order to compute $E^2_{1,1}$ we again use
Poincar{\'e}-duality and find
\[E^2_{1,1} \cong\rmH_1(G',\vrh')
\cong \rmH^1(G',\vrh').\]
\nid So we need to determine the
derivations of $\vrh'$. We make an Ansatz
\[ d(m)=\mu={r\choose s},\;\;\;\;d(\ell)=\lam={t\choose u}.\]
\nid This defines a derivation if and only if $0=d(1)=d(m\ell
m^{-1}\ell^{-1})$ which by virtue of the relations $d(gh)=dg+g\cdot
dh$ and $d(g^{-1})=-g^{-1}\cdot dg$
transforms by expansion into the equation
\begin{eqnarray*}
(\vrh'(m)-\Id)\lam&=&(\vrh'(\ell)-\Id)\mu\\[0.8ex]
\left(\begin{array}{rr}0&1\\0&0\end{array}\right){r\choose s}&=&
\left(\begin{array}{rr}-2&6\\0&-2\end{array}\right){t\choose u}.
\end{eqnarray*}
\nid The most general solution is given by choosing
$r,t\in\ZZ$ arbitrary and setting $s=0,\; u=-2t$. Therefore
$\Der(\vrh')\cong\ZZ^2$. The corresponding derivation $d_{r,t}$
given by $\mu=(r,0)$ and
$\lam=(t,-2r)$ is principal if and only if there is
$\kap\in\ZZ^2$ so that
\[\mu=(\vrh'(m)-\Id)\kap\;\;\;\mathrm{and}\;\;\;
\lam=(\vrh'(\ell)-\Id)\kap.\]
\nid The unique solution to these equations is
\[\kap={-\frac{t}{2}-3r\choose r}\]
\nid which is integral only if $t\equiv0\MOD 2$, hence the principal
derivations correspond to the submodule $\ZZ\oplus 2\ZZ\subset\ZZ^2$
and so
\[\rmH_1(G',\vrh')
\cong\rmH^1(G',\vrh')
\cong\left(\frac{\ZZ^2}{2\ZZ\oplus\ZZ}
  \right)  \cong
\ZZ/2.\]

\nid Thus we find the term
\[E^2_{1,1}\cong \ZZ/2\]
\nid which completes the computation of the last missing term in the
spectral sequence. Again
the sequence degenerates because $E^2_{2,1}=\{0\}$, and the other arrow
$\d_2\colon E^2_{2,0}\to E^2_{0,1}$ vanishes because
\[\rmH_2(B_0\cap
B_1,\ZZ)\stackrel{\sig_*}{\longrightarrow}
E^2_{2,0}\stackrel{\d_2}{\longrightarrow}
E^2_{0,1}\]
\nid is an exact sequence. Here $\sig$ is a section with
$\sig_*=(f_*|E^2_{2,0})^{-1}$, hence $\sig_*$ is surjective and
$\d_2=0$. The homology is thus given
by the total homology of the $E^2$-term which is as described in the
proposition. \epf

\th{TS3}{The integral homology of $Y$ is as follows:
\begin{center}
\begin{tabular}{|c|cccccc|}
\hline
$p$&$0$&$1$&$2$&$3$&$4$&$5$\\
\hline
$\rmH_p(Y,\ZZ)$&$\ZZ$&$0$&$\ZZ$&$\ZZ$&$0$&$\ZZ$\\
\hline
\end{tabular}
\end{center}}\eth

\pf We apply the Mayer-Vietoris sequence
\begin{equation}\label{MVVV}
\cdots\to
\rmH_{p+1}(Y,\ZZ)\to\rmH_p(Y_0\cap Y_1,\ZZ)\stackrel{\psi_p}{\to}
\rmH_p(Y_0,\ZZ)\oplus
\rmH_p(Y_1,\ZZ)\to \rmH_p(Y,\ZZ)\to\cdots\end{equation}
\nid to $Y=Y_0\cup Y_1$.
Let 
\[[\mu],[\lam]\in\pi_1(B_0\cap B_1)\cong H_1(B_0\cap B_1,\ZZ)\]
be the classes of a meridian and
a longitude, respectively, let
$F$ be a fibre and $\sig$ a section of $Y_0\to B_0$. Now obviously
$H_0(Y,\ZZ)=\ZZ$. We also know already that
$H_1(Y,\ZZ)=\{0\}$ because \reft{Psc} states that $Y$ is simply
connected. But observe that $H_1(Y_0\cap Y_1,\ZZ)_{\mathrm{free}}\cong\ZZ^2$
is generated by $\sig_*[\mu]$ and $\sig_*[\lam]$, while
$H_1(Y_0,\ZZ)\cong\ZZ(\sig_*[\mu])$  and
$H_1(Y_1,\ZZ)_{\mathrm{free}}\cong\ZZ(\sig_*[\lam])$. Then 
\[\psi_1
( \sig_*[\mu]) = (\sig_*[\mu],0)\;\; \mathrm{and}\;\; \psi_1
( \sig_*[\lam]) = (0,-\sig_*[\lam])\in \rmH_1(Y_0,\ZZ)\oplus H_1(Y_1,\ZZ).\]
This implies that $\psi_1$ induces an isomorphism of the free parts
\begin{equation}\label{eettt}\rmH_1(Y_0\cap Y_1,\ZZ)_{\mathrm{free}}
\stackrel{\cong}{\longrightarrow} \rmH_1(Y_0,\ZZ)\oplus
\rmH_1(Y_1,\ZZ)_{\mathrm{free}}.\end{equation}
\nid Let $[\tau_1]\in H_1(Y_1,\ZZ)$ and $[\tau_2]\in H_1(Y_0\cap
Y_1,\ZZ)$ be the unique torsion elements.
Since $\psi_1$ has to be onto in order for $H_1(Y,\ZZ)$ to vanish, we
conclude $[\tau_1]\in\im(\psi_1)$. Given the isomorphism (\ref{eettt})
this is only possible if $\psi_1([\tau_2])=[\tau_1]$
for the torsion class, and
thus $\psi_1$ is an isomorphism of $\ZZ$-modules. So
$\coker(\psi_1)=\ker(\psi_1)=\{0\}$.
Therefore $H_2(Y,\ZZ)\cong\coker(\psi_2)$. Note that
\[\rmH_2(Y_0,\ZZ)_{\mathrm{free}} \cong \rmH_2(Y_1,\ZZ)\cong\ZZ[F]\]
\nid is generated by the class of a fibre, while for a generator $T\in
H_2(B_0\cap B_1,\ZZ)$
\[\rmH_2(Y_0\cap Y_1,\ZZ)_{\mathrm{free}}  \cong
\ZZ[F] \oplus \ZZ(\sig_*[T]).\]
\nid Thus $\psi_2([F])=([F],-[F])$ while
$\psi_2(\sig_*[T])=0$. Now denote the degree-$2$
torsion classes by $[\tau_3]\in H_2(Y_0,\ZZ)$ and $[\tau_4]\in
H_2(Y_0\cap Y_1,\ZZ)$. Suppose $\psi_2([\tau_4])=[\tau_3]$. In this case
we get
\begin{eqnarray}
\label{ee12}\ker(\psi_2)&=&\ZZ([\sig_*T])\cong\ZZ,\\
\label{ee13}\coker(\psi_2)&=&\frac{{\ZZ[F] \oplus
  \ZZ[F]\oplus \ZZ/2[\tau_3]}}{{\{(a[F],-a[F])\;|
  \;a\in\ZZ\}\oplus\ZZ/2[\tau_3]}} \cong \frac{\ZZ^2}{\diag(\ZZ^2)}\cong\ZZ
\end{eqnarray}
\nid and therefore
\[\rmH_2(Y,\ZZ)\cong\coker(\psi_2)\cong\ZZ[T]\cong\ZZ.\]
\nid Otherwise if $\psi_2([\tau_4])\neq\tau_3$ then
$\psi_2([\tau_4])=0$ since it is a
torsion class. But then we obtain
$H_2(Y,\ZZ)\cong\ZZ\oplus\ZZ/2$. This is a contradiction since the
torsion subgroup $H_2(Y,\ZZ)_{\mathrm{tors}}\cong \half T\oplus \half
T$ of a simply connected spin five-manifold
is always decomposable as a direct sum of two isomorphic torsion
modules (see \cite{Smale}). So we can rule out the case
$[\tau_4]\in\ker(\psi_2)$ and
$H_2(Y,\ZZ)\cong\ZZ$. For the remaining homology groups $H_p(Y,\ZZ)$
with $p>2$ we only give a sketch since Smale's \reft{TClass} implies
that they are already uniquely determined by $H_2(Y,\ZZ)$.
For $Y=Y_0,Y_1,Y_0\cap Y_1$ there is an isomorphism
\[\tau_{[F]}\colon \rmH_1(Y,\ZZ)_{\mathrm{free}}\longrightarrow
\rmH_3(Y,\ZZ).\]
\nid This isomorphism sends the classes of $m$ and $\ell$ to
$[f^{-1}(m)]$ and $[f^{-1}(\ell)]$. Then
\[\psi_3\circ \tau_{[F]}=(\tau_{[F]},\tau_{[F]})\circ\psi_1,\]
\nid hence $\psi_3$ is an isomorphism as is $\psi_1$. Thus
\[\rmH_3(Y,\ZZ)=\ker(\psi_2)\cong\ZZ\]
\nid by (\ref{ee13}). For $i>3$, we have $\psi_i=0$ trivially, and
\[\rmH_4(Y,\ZZ)=\coker(\psi_4)=\frac{\{0\}}{\{0\}}=\{0\}\]
\nid since $H_4(Y_i,\ZZ)=\{0\}$ for
$i=0,1$. Finally
\[\rmH_5(Y,\ZZ)=\ker(\psi_4)=\rmH_4(Y_0\cap Y_1,\ZZ)\cong \ZZ,\]
\nid which is also clear otherwise since $Y$ is compact, orientable
and $\dim(Y)=5$.\epf

\cor{Chs3}{$Y$ is an (integral) homology $S^2\times S^3$.}\ecor

We are now in a position to show that $Y$ is diffeomorphic to
$S^2\times S^3$. This uses the following theorem of Smale:

\th{TClass}{{\rm{(Smale \cite{Smale}) }} There is a bijective correspondence
  between the category of smooth, closed, orientable, simply connected,
  $5$-dimensional spin manifolds and the category of finitely
  generated Abelian groups. The correspondence is realized by sending
  a manifold $Y$ in this category to the group $F\oplus\half T$, where
  $F\oplus T=\rmH_2(Y,\ZZ)$ is a decomposition of the second integral
  homology of $Y$ into its free and its torsion subgroups, and
  $T=\half T \oplus\half T$ is a direct sum decomposition of the
  torsion subgroup.}\eth

This is an analogue for five-manifolds of Wall's classification
theorem in dimension six, however, 
it is easier to apply since we do not need to know
the ring-structure on homology, and it is also more general in that it
includes manifolds with
torsion in homology.

\th{TDif}{$Y$ is diffeomorphic to $S^2\times S^3$.}\eth
\pf The five-manifold $Y$ is smooth and compact
by \refp{Psmo} and simply connected by
\reft{Psc}, it is orientable by \refp{Pori}.
$Y$ is torsion free and has the right Betti-numbers
by \reft{pcht23}. $Y$ is
spin from \reft{TSpin}.
Hence the claim follows by applying \reft{TClass}.\epf

This gives us torus fibrations on well-known manifolds. We
formulate this in the following theorem.

\th{TT3f}{There are $T^3$-fibrations on $S^1\times
  S^2\times S^3$ and $S^3\times S^3$ whose singularities are locally
  equivalent to (\ref{Diag}).}\eth

\pf For the first we simply take the product $S^1\times Y$ and extend
the fibration map $f\colon Y\to S^3$ trivially, for the second we take
the Hopf map $H\colon S^3\to S^2$ and let $\tilde{f}\colon S^3\times S^3\to
S^3$ be given by $\tilde{f}=f\circ(H,\id)$. Note that this latter fibration
has no section, while the first one obviously has one.\epf

\ssec{S4SS2}{General torus knots}

We give a brief account of the general torus knot $\grt(p,q)$ for arbitrary 
$p,q\in\NN$.
First observe that we can easily modify the argument given in
\refss{S4SS1} to obtain a smooth compactification of a $T^2$-bundle
over $S^3\setminus \grt(p,q)$ for any pair of coprime integers
$(p,q)\in\NN^2$
such that $2|p$ and $3|q$. As in (\ref{eqeqeqeqeq}), we define
$Y_{p,q}\subset \CC P^2\times S^3$ by the (single-valued) equation
\[4x^3-g_2^{q/3}xz^2-g_3^{p/2}z^3-y^2z=0\]

\nid and let again the fibration $f\colon Y_{p,q}\to S^3$ be
given by $f={\mathrm{pr}}_2|Y_{p,q}$ where
${\mathrm{pr}}_2\colon \CC P^2\times S^3\to S^3$ is projection onto
the second
factor. As before, a local calculation shows $Y_{p,q}$ 
to be a smooth compact submanifold
of $\CC P^2\times S^3$.
For generic $(g_2,g_3)\in S^3$ the fibre $f^{-1}(g_2,g_3)$ is
an elliptic curve in $\CC P^2$ determined
by the Weierstra{\ss} equation with coefficients $g_2^{q/3}$
and $g_3^{p/2}$. The fibre is singular if and only if the
discriminant vanishes, that is
\[\Del=27\left(g_3^{p/2}\right)^2-\left( g_2^{q/3}
\right)^3 =27 g_3^p-g_2^q=0.\]

\nid From \cite{Mi2} the discriminant locus is a $(p,q)$-torus knot.
The argument for the vanishing of the second Stiefel-Whitney
class in \refss{S4SS1} carries over directly to $Y_{p,q}$. 
From \reft{Psc} we deduce that
$Y_{p,q}$ is simply connected. To summarize:

\th{TT43pq}{$Y_{p,q}\subset\CC P^2\times S^3$ is a smooth, compact,
  orientable
  submanifold with $\pi_1(Y_{p,q})=\{1\}$ and $w_2(Y_{p,q})=0$. There
  is a $T^2$-fibration $f\colon Y_{p,q}\to S^3$ with singularities of
  type (\ref{Diag}) which
  degenerates over a $(p,q)$-torus knot $\grt(p,q)\subset S^3$.}\eth

% ^ 

\subsubsection{$M$-representations}

Using the techniques developed in \refs{S3} one can show that  a torus knot 
$\grt(p,q)\subset S^3$ admits a non-trivial 
$M$-representation (that is, one whose image is $\ncong\ZZ$) if and
only if 
$2|p$ and $3|q$ (after exchanging $p$ and $q$ if necessary). In this case,
there is precisely one $M$-representation whose image in $\SL(3,\ZZ)$ is 
isomorphic to $\SL(2,\ZZ)$. This is the monodromy of the $T^2$-fibration 
$Y_{p,q}\to S^3$ described above.
The case $p=2, \;q=3$ was dealt with in \reft{P2345610}. We give a short
description of the other $\SL(2,\ZZ)$-valued $M$-representations which occur 
as monodromy representations in the $T^2$-fibrations $Y_{p,q}\to S^3$.
We refer to \refas{A2} for more details and a complete classification of all
$M$-representations.\\

\pid \reft{T43}
and its proof show that the group
$G=G(\grt({4,3}))$ of the 
$(4,3)$-torus knot $\grt(4,3)$ has
a unique non-trivial $M$-re\-pre\-sen\-ta\-tion 
$\vrh=\vrh_0\colon G\to \SL(2,\ZZ)\subset \SL(3,\ZZ)$ which
can be given, with respect to a presentation $G=\<x,y\;|\;x^4=y^3\>$, by
\begin{equation}\label{eqt43}
\vrh(x)=\left(\begin{array}{rrr}0&1&0\\-1&0&0\\0&0&1\end{array}\right),
   \;\;\;\;
\vrh(y)=\left(\begin{array}{rrr}0&1&0\\-1&-1&0\\0&0&1\end{array}\right),
\end{equation}

\nid see \refas{A2}. 
Note that $\vrh(x)^4=\vrh(y)^3=\Id$ 
and the meridian $m=x^{-1}y$ is mapped to the
standard monodromy matrix $A$ given in (\ref{eeaa}). \\

\pid Now suppose we have $p,q\in\NN$ with $2|p$, $3|q$, $4\nmid p$ and 
$\gcd(p,q)=1$. 
Let $G(\grk)=\<x,y\,|\,x^p=y^q\>$ and
$G{(2,3)}=\<\bar{x},\bar{y}\,|\, \bar{x}^2=\bar{y}^3\>$. Let
$\bar{\vrh}\colon G{(2,3)}\to\SL(2,\ZZ)$ the unique non-trivial
$M$-re\-pre\-sen\-ta\-tion of $G{(2,3)}$. Then 
\[\vrh(x)=\bar{\vrh}(\bar{x})^{\pm1}\;\;\;\;\mathrm{ and }\;\;\;\;
\vrh(y)=\bar{\vrh}(\bar{y})^{\pm1}.\]
Also note that the image of a
longitude $\ell$ under our re\-pre\-sen\-ta\-tion is 
\[\vrh(\ell)=\left(\begin{array}{ccc}-1&pq&0\\ 0&-1&0\\
0&0&1\end{array}\right)=
\bar{\vrh}(\bar{\ell})^{pq/6}.\]

If  $p,q\in\NN$ with $4|p$ and $3|q$ the unique $M$-representation with 
image $\SL(2,\ZZ)$ can be derived similarly from the $\SL(2,\ZZ)$-valued
$M$-representation (\ref{eqt43}) of $G(\grt(4,3))$.

% ^ 
\subsubsection{$\grt(4,3)$  and 
$(S^3\times S^3)\#(S^3\times S^3)\#(S^4\times S^2)$}

We now study the $T^2$-bundle
$f\colon Y=Y_{4,3}\to S^3$ associated to the $\SL(2,\ZZ)$-valued
$M$-re\-pre\-sen\-ta\-tion of $G(\grt(4,3))$. As before, let
\[\begin{array}{lclclcl}
B_0&=&S^3\setminus \grt{(4,3)},&&Y_0=f^{-1}(B_0),\\
B_1&=&\calN(\grt{(4,3)}),&&Y_1=f^{-1}(B_1).
\end{array}\]
\nid and denote by $M$ the $\ZZ G$-module given by the monodromy
re\-pre\-sen\-ta\-tion $\vrh$. Since the computations are somewhat lengthy 
and all the necessary techniques were already illustrated in the last 
example, we will only give the results.

\th{Pcoh43}{The integral homology of $Y_0$, $Y_1$ and $Y_0\cap Y_1$
  is given by the following table:
\begin{center}
\begin{tabular}{|c||c|c|c|}
\hline
$p$&$\rmH_p(Y_0)$&$\rmH_p(Y_1)$&$\rmH_p(Y_0\cap Y_1)$\\
\hline\hline
$0$&$\ZZ$           &$\ZZ$           &$\ZZ$\\
$1$&$\ZZ$           &$\ZZ^2$         &$\ZZ^3$\\
$2$&$\ZZ^3$         &$\ZZ^2$         &$\ZZ^4$\\
$3$&$\ZZ^3$         &$\ZZ$           &$\ZZ^3$\\
$4$&$0$             &$0$             &$\ZZ$\\
$5$&$0$             &$0$             &$0$\\
\hline
\end{tabular}
\end{center}}\eth

By applying the Mayer-Vietoris sequence one finds 
the homology of $Y=Y_0\cup Y_1$.

\th{Tc43}{The integral homology of $Y=Y_0\cup Y_1$ is as follows:
\begin{center}
\begin{tabular}{|c|cccccc|}
\hline
$p$&$0$&$1$&$2$&$3$&$4$&$5$\\
\hline
$\rmH_p(Y,\ZZ)$&$\ZZ$&$0$&$\ZZ^2$&$\ZZ^2$&$0$&$\ZZ$\\
\hline
\end{tabular}
\end{center}}\eth

The following corollary follows from Smale's \reft{TClass} in the
same way as the 
corresponding statements in \refss{S4SS1}.

\cor{Pcss}{$Y$ is diffeomorphic to $(S^2\times S^3)\#(S^2\times
  S^3)$.}
\ecor

\pf $Y$ is closed, spin, and simply connected by \reft{TT43pq}. The
description of the integral homology in \reft{Tc43} together with
Smale's \reft{TClass} yields the result.\epf

\nid As 
before we can construct $\TT^3$-fibrations $X\to Y\to S^3$ where 
$X\to Y$ is a $\TT^1$-bundle. If we chose the trivial $\TT^1$-bundle
then
\[X\cong S^1\times(S^2\times S^3)\#(S^2\times S^3),\]
and if $X\to Y$ is a $\TT^1$-bundle with primitive Chern-class in 
$H^2((S^2\times S^3)\#(S^2\times S^3),\ZZ)\cong\ZZ^2$ then
we obtain the space
\[X\cong (S^3\times S^3)\#(S^3\times S^3)\#(S^4\times S^2).\]
This uses again the Leray spectral sequence.

\subsubsection{The $(p,q)$-torus knot}
We argue that we do not get any new
manifolds by using other torus knots.

\th{Tothers}{Let $p,q\in\NN$ be coprime positive integers with $2|p$
  and $3|q$, let $\grk=\grt(p,q)\subset S^3$ be the $(p,q)$-torus knot
  and let $f\colon Y_{p,q}\to S^3$ be the $T^2$-bundle with
  section associated to the unique $M$-re\-pre\-sen\-ta\-tion
  $\vrh=\vrh_0\colon 
  G(\grk)\to \SL(2,\ZZ)\subset\SL(3,\ZZ)$. 
  Then either $4\nmid p$ and $Y_{p,q}\cong S^2\times
  S^3$ or $4|p$ and $Y_{p,q}\cong (S^2\times S^3)\#(S^2\times S^3)$.}\eth

{\it Sketch of proof. } One can then check that all calculations 
for $\grt(2,3)$ and $\grt(4,3)$ remain valid. \epf

\cor{Pinfi}{There are infinitely many pairwise inequivalent 
  $T^3$-fibrations with singularities of type (\ref{Diag})
\begin{eqnarray}
S^3\times S^3&\longrightarrow&S^3\label{e111}\\
S^1\times S^2\times S^3&\longrightarrow&S^3\label{e222}\\ 
S^1\times((S^2\times S^3)\#(S^2\times S^3))&\longrightarrow&S^3\label{e333}\\
(S^3\times S^3)\#(S^3\times S^3)\#(S^4\times S^2)&
\longrightarrow&S^3.\label{e444}
\end{eqnarray}
  The discriminant locus can be any torus knot $\grt(2p',3q')$, with $p'$ odd
  in (\ref{e111}), (\ref{e222}) and $p'$ even in (\ref{e333}), (\ref{e444}).
  More generally such fibrations exist on every $\TT^1$-bundle 
  over $S^2\times S^3$ or
  $(S^2\times S^3)\#(S^2\times S^3)$.}\ecor

\nid 

\appendix

\sec{A2}{ Computations with knot groups}

\th{T43}{Let $\grk\subset S^3$ be the $(4,3)$-torus knot and
  $G=\pi_1(S^3\setminus\grk)$ its group. Then the
  $M$-re\-pre\-sen\-ta\-tions of 
  $G$ are given by the trivial $M$-re\-pre\-sen\-ta\-tion
  $\varrho_{\mathrm{trivial}}$ and an infinite family $\varrho_k$ 
  of non-Abelian $M$-re\-pre\-sen\-ta\-tions parametrised by
  $k\in\ZZ$, with $\varrho_{-k}=(\varrho_k^{-1})^{T}=\varrho_k^*$. }\eth

\begin{pspicture}(1,0)(16,6)
\pscustom[linewidth=1.5pt]{\pscurve(8.2,2.0)(7.5,1.8)(7.0,1.9)(6.5,3.0)
  (6.5,4.5)(7.0,5.4)(8.0,6.0)(9.0,5.8)(9.8,4.8)(9.6,4.1)}
\pscustom[linewidth=1.5pt]{\pscurve(9.4,3.6)(9.3,3.1)(10.0,2.6)}
\pscustom[linewidth=1.5pt]{\pscurve(10.4,2.2)(10.6,1.7)(9.5,1.1)(8.5,2.0)
  (8.0,4.0)(8.2,4.4)}
\pscustom[linewidth=1.5pt]{\pscurve(8.4,4.9)(8.6,5.4)(8.8,5.6)}
\pscustom[linewidth=1.5pt]{\pscurve(9.3,5.9)(10.0,6.0)(10.5,5.8)(10.8,5.4)
  (10.5,4.2)(9.5,3.9)(8.5,3.4)(8.4,3.3)}
\pscustom[linewidth=1.5pt]{\pscurve(7.8,3.2)(7.2,3.8)(7.5,4.5)(9.5,4.6)}
\pscustom[linewidth=1.5pt]{\pscurve(10.0,4.6)(10.2,4.7)(10.4,5.0)(10.6,5.1)}
\pscustom[linewidth=1.5pt]{\pscurve(11.1,5.2)(11.4,5.0)(11.3,3.0)(8.8,2.1)}
\psline[linewidth=0.8pt,arrows=->](7.7,2.4)(8.1,2.4)
\pscurve[linewidth=0.8pt,arrows=->](8.5,2.4)(8.7,2.5)(8.5,2.7)(7.7,2.6)
\rput(7.4,2.4){$g_1$}
\psline[linewidth=0.8pt,arrows=->](8.0,5.2)(8.3,5.2)
\pscurve[linewidth=0.8pt,arrows=->](8.7,5.2)(9.0,5.3)(8.8,5.5)(8.0,5.4)
\rput(7.6,5.3){$g_8$}
\psline[linewidth=0.8pt,arrows=->](9.2,3.0)(8.8,2.9)
\pscurve[linewidth=0.8pt,arrows=->](8.8,3.1)(9.6,3.4)(9.8,3.2)(9.6,3.1)
\rput(10.1,3.2){$g_6$}
\psline[linewidth=0.8pt,arrows=->](11.3,3.5)(10.9,3.5)
\pscurve[linewidth=0.8pt,arrows=->](10.9,3.7)(11.7,3.8)(11.9,3.6)(11.7,3.5)
\rput(12.2,3.6){$g_3$}
\psline[linewidth=0.8pt,arrows=->](5.9,3.6)(6.3,3.6)
\pscurve[linewidth=0.8pt,arrows=->](6.7,3.6)(6.9,3.7)(6.7,3.9)(5.9,3.8)
\rput(5.4,3.4){$g_2$}
\psline[linewidth=0.8pt,arrows=->](11.4,5.5)(11.0,5.5)
\pscurve[linewidth=0.8pt,arrows=->](10.6,5.5)(10.4,5.7)(10.6,5.9)(11.4,5.7)
\rput(11.8,5.6){$g_5$}
\psline[linewidth=0.8pt,arrows=->](10.3,5.5)(10.4,5.1)
\pscurve[linewidth=0.8pt,arrows=->](10.4,4.7)(10.3,4.5)(10.2,4.6)(10.1,5.5)
\rput(9.9,5.7){$g_7$}
\psline[linewidth=0.8pt,arrows=->](7.0,4.8)(7.3,4.5)
\pscurve[linewidth=0.8pt,arrows=->](7.5,4.3)(7.5,4.2)(7.4,4.1)(6.9,4.6)
\rput(7.6,3.9){$g_4$}
\pscurve[linewidth=0.8pt,linestyle=dashed,arrows=->](12.8,2.3)(11.8,2.7)
(11.4,2.8)
\pscurve[linewidth=0.8pt,linestyle=dashed,arrows=->]
(10.9,2.9)(10.3,2.9)(10.1,2.8)  
\pscurve[linewidth=0.8pt,linestyle=dashed,arrows=->](9.6,2.6)(9.4,2.6)
(9.0,3.5) 
\pscurve[linewidth=0.8pt,linestyle=dashed,arrows=->](9.0,3.9)(9.2,4.2)
(10.1,4.3)(10.8,3.4)(12.8,2.5)
\rput(13.2,2.5){$x$}
\setcounter{psfig}{1}
\rput(8.0,0.5){{\bf Figure 2: }
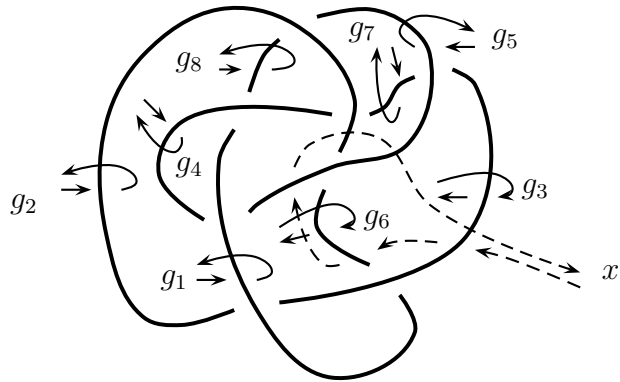
{\it{The $(4,3)$-torus knot $\grt{(4,3)}$ (or $8_{19}$).}}}
\end{pspicture}

We do not give a full proof which is very similar to that of \reft{P2345610}, 
but we do describe explicitly these $M$-representations.
From the above knot-diagram we obtain the following
Wirtinger-presentation of $G$:\\%[-10.2mm]
\nopagebreak
\setlength{\fboxsep}{\tabcolsep}
\setlength{\tabcolsep}{-3mm}
\makebox[\width][r]{
\begin{tabular}{p{0.5\textwidth}p{0.5\textwidth}}
\begin{eqnarray}
\label{G1} g_1^{-1}g_3g_1&=&g_2 \\
\label{G2} g_3^{-1}g_6g_3&=&g_1 
\end{eqnarray}
&\hfill
\begin{eqnarray}
\label{G5} g_2^{-1}g_8g_2&=&g_5\\
\label{G6} g_4^{-1}g_1g_4&=&g_8
\end{eqnarray}
\end{tabular}}
\setlength{\tabcolsep}{\fboxsep}\\

\vspace*{-14mm}
\setlength{\fboxsep}{\tabcolsep}
\setlength{\tabcolsep}{-3mm}
\makebox[\width][r]{
\begin{tabular}{p{0.5\textwidth}p{0.5\textwidth}}
\begin{eqnarray}
\label{G3} g_5^{-1}g_2g_5&=&g_6 \\
\label{G4} g_2^{-1}g_4g_2&=&g_7 
\end{eqnarray}
&\hfill
\begin{eqnarray}
\label{G7} g_1^{-1}g_5g_1&=&g_4\\
\label{G8} g_5^{-1}g_7g_5&=&g_3
\end{eqnarray}
\end{tabular}}
\setlength{\tabcolsep}{\fboxsep}
\vspace*{-4.8mm}

Making an Ansatz as in the case of a trefoil 
we get the family $(\varrho_k)_{k\in\ZZ}$ of $M$-representations, 
determined by sending
\[g_1\longmapsto
\left(\begin{array}{ccc}1&1&0\\0&1&0\\0&0&1\end{array}\right),\;\;
g_5\longmapsto\label{mmr}
\left(\begin{array}{ccc}1&0&0\\-1&1&0\\0&0&1\end{array}\right)\]

\nid and sending the remaining generator
\[g_3\longmapsto
\left(\begin{array}{ccc}0&1&k\\-1&2&k\\0&0&1\end{array}\right)\;\;
\mathrm{if\;}k<0,\mathrm{\;\;\;and\;\;\;}
g_3\longmapsto
\left(\begin{array}{ccc}0&1&0\\-1&2&0\\k&-k&1\end{array}\right)\;\;
\mathrm{if\;}k\geq0.\]
\epf

\nid For $k=0$ we get an $M$-representation with image $\SL(2,\ZZ)\subset
\SL(3,\ZZ)$. This is the unique non-trivial $M$-representation of 
$G(\grt(4,3))$ with this property and was used in \refss{S4SS2}.

\pid Similar computations show that the possible $M$-representations of a 
general torus knot can be obtained from those of the $(2,3)$ and $(4,3)$ 
torus knots. We do not carry this out here but only state the result.

\th{PT}{The group $G{(p,q)}$\label{gpq} 
  of the $(p,q)$-torus knot $\grt{(p,q)}$ has non-trivial
  $M$-re\-pre\-sen\-ta\-tions if and only if
  (after perhaps exchanging $p$ and $q$)
  $p=2^mp'$ and $q=3^kq'$ with $\gcd(p',6)=\gcd(q',6)=\gcd(p',q')=1$ and
  $m,k>0$. In case $m=1$ there is a bijection
  $\MRep(G{(p,q)})\to \MRep(G{(2,3)})$ between the
  $M$-re\-pre\-sen\-ta\-tions of 
  $G{(p,q)}$ and $G{(2,3)}$. In case $m>1$ there is a bijection
  $\MRep(G{(p,q)})\to \MRep(G{(4,3)})$.}\eth

If $G=\<x,y\;|\;x^2=y^3\>$ is the trefoil group and $\bar{G}=\<\bar{x}, 
\bar{y}\;|\;\bar{x}^{2p'}=\bar{y}^{3q'}\>$ the group of the $(2p',3q')$ 
torus knot, with $p',q'$ odd and relative prime,
this correspondence is given by sending an $M$-representation $\vrh$ of $G$ to
the $M$-representation $\bar{\vrh}$ of $\bar{G}$ which is defined by
\begin{equation}\label{rewww}
\bar{\vrh}(\bar{x})=\vrh(x)^{\pm1},\;\;\;\;
\bar{\vrh}(\bar{y})=\vrh(y)^{\pm1},
\end{equation}
for a unique choice of exponent $\pm1$ which is determined by $p'$ and $q'$.
This is well-defined because $\vrh(x)^2=\vrh(y)^3=-\Id$ and $p',q'$ are both 
odd, hence (\ref{rewww}) does define a homomorphism, and there is a unique 
choice of exponents $\pm1$ in (\ref{rewww}) to make it an $M$-homomorphism. 
The 
case of $\grt(4p',3q')$ with $\gcd(2p',q')=1$ is analogous.

\bibliographystyle{alpha}

\end{document}